\documentclass{amsart}

\usepackage{amsmath}

\usepackage{amssymb}
\theoremstyle{plain}   
\begingroup 

\newtheorem{theorem}{Theorem}[section]   
\newtheorem{lemma}[theorem]{Lemma}         
\newtheorem{proposition}[theorem]{Proposition}  
\endgroup

\theoremstyle{definition}
\newtheorem{definition}[theorem]{Definition}   

\theoremstyle{remark}
\newtheorem{remark}[theorem]{Remark}        
\newtheorem{example}[theorem]{Example}        

\numberwithin{equation}{section}

\newcommand{\R}{{\mathbb R}}
\newcommand{\N}{{\mathbb N}}
\newcommand{\Z}{\mathbb Z}

\newcommand{\e}{\varepsilon}

\newcommand{\ep}{\varepsilon}

\newcommand{\vf}{\varphi}

\newcommand{\interior}{\operatorname{int}}
\newcommand{\inte}{\operatorname{int}}

\newcommand{\Hau}{{\mathcal{H}}}

\newcommand{\mcA}{\mathcal{A}}

\newcommand{\p}{{\mathcal{P}}}
\newcommand{\s}{{\mathcal{S}}}
\newcommand{\Ka}{{\mathcal{K}}}

\newcommand{\wtilde}{\widetilde}

\newcommand{\mcI}{{\mathcal{I}}}
\newcommand{\mcJ}{{\mathcal{J}}}

\newcommand{\mcP}{{\mathcal{P}}}

\usepackage[dvips]{color}

\begin{document}

\title[Curves in Banach spaces]{Curves in Banach spaces which allow a $C^2$ parametrization}

\author{Jakub Duda}
\author{Lud\v ek Zaj\'\i\v{c}ek}



\email{jakub.duda@gmail.com}
\email{zajicek@karlin.mff.cuni.cz}

\address{
\v{C}EZ, a. s.,
Duhov\'a 2/1444,
140 53 Praha 4,
Czech Republic} 

\address{
Charles University,
Faculty of Mathematics and Physics,
Sokolovsk\'a 83,
186 75 Praha 8-Karl\'\i n,
Czech Republic}

\begin{abstract}
We give a complete characterization of those $f: [0,1] \to X$ (where $X$
is a Banach space which admits an equivalent Fr\'echet smooth norm) which allow an equivalent $C^2$ parametrization.
For $X=\R$, a characterization  is well-known.
However, even in the case $X=\R^2$,
several quite new ideas are needed. Moreover, the very close case of parametrizations with a bounded second derivative
 is solved.
\end{abstract}

\maketitle

\section{Introduction} 

Let $X$ be a (real) Banach space, and a curve $f:[a,b]\to  X$ be
given. Several authors (e.g.  Ward, Zahorski, Choquet, Tolstov)
 investigated  conditions under which $f$~allows an equivalent
parametrization which is ``smooth of the first order'' 
(e.g., it is differentiable, boundedly differentiable, or continuously differentiable). 
Their results, which deal with
the case $X=\R^n$,  were generalized (using more or
less difficult modifications of known methods) to the case of an arbitrary $X$
in~\cite{DZ} and~\cite{D}, where information concerning the history of the 
``first order case'' can be found.
\par
The case of ``higher order smooth'' parametrizations for
$X=\R$ was settled independently by Laczkovich and Preiss \cite{LP}
and Lebedev~\cite{Leb}. Both papers contain (formally slightly
different, see Section 7 below) characterizations of those $f: [0,1] \to \R$ which
allow an equivalent $C^n$ ($n \in \N)$ parametrization or a $C^{n,\alpha}$ ($0 < \alpha \leq 1$) parametrization
 (i.e., a parametrization, whose $n$-th derivative is $\alpha$-H\" older). The case of $f: [0,1] \to \R$
 which allow an $n$-times  differentiable parametrization was settled in \cite{Duge}.
 \par
The problem of ``higher order smooth'' parametrizations
in the vector case (even for $X=\R^2$) is essentially more difficult
than in the case $X=\R$.  In the present article we characterize those $f: [0,1] \to X$ which
allow an equivalent $C^2$ parametrization if $X$  admits an equivalent  Fr\' echet smooth norm. We were not able to solve the vector-valued problem
of $C^n$ parametrization for $n\geq 3$ (even for $X=\R^2$). 
It seems to us that if there is a satisfactory solution of this problem, 
 then it requires some new ideas.
\par
An earlier (unpublished) preprint [DZ2] treated the case of $C^2$ parametrizations together with the case of parametrizations
with bounded convexity. 
In the present article, the more interesting $C^2$ case is treated separately and with a simplified proof of
 a basic lemma (Lemma \ref{Wfin2} below).  Moreover, we also treat a very similar case of curves
allowing a $D^{2,\infty}$ parametrization (i.e., parametrization with bounded second derivative), since we observed
 that the main lemmas (Lemma \ref{Wfin2} and Lemma \ref{hlavkon2}) can be proved for both cases simultaneously. The
  only new idea in the  $D^{2,\infty}$ case is the definition of the notion of  an $[f,\delta,K]$-partition
   (together with Lemma \ref{konstrpartlem} and Lemma \ref{nove}) which is more complicated than that of 
    an $(f,\delta,K)$-partition (which is sufficient for the $C^2$ case). Let us remark that the $C^2$ and
     $D^{2,\infty}$ cases are equivalent for $X=\R$ (see \cite[Remark 3.7]{LP} or Section 7 below), but are 
      not equivalent for $X=\R^2$ (see Example \ref{dmd} below).
\par
An essential modification of our methods gives a full characterization of vector-valued curves allowing 
$C^{1,\alpha}$ ($0 < \alpha \leq 1$) parametrizations \cite{DZ3}. Using another method (which follows
 \cite{Leb} and \cite{LP}), the case of of vector-valued curves allowing  a twice differentiable parametrizations
  was settled in \cite{Duse}. 
\par
The structure of the paper is the following. In Section 2 we recall some known notions and their 
well-known properties.
 In Section 3 we introduce the main special notions (the variation $W^{\delta}(f,G)$ and some types of
  generalized partitions of open sets) and prove a number of lemmas. 
The core of the article are Lemma \ref{Wfin2} and Lemma~\ref{hlavkon2}.
 Section 4 contains main results on $C^2$ parametrizations. A consequence,
 which might  be interesting also from the point of view of differential geometry, is presented in Proposition \ref{klas}.
In Section 5  the results concerning $D^{2,\infty}$ parametrizations are presented (the proofs which are almost literally same as those in Section 4 are omitted). Section 6 contains  several examples, which show applicability of our results.
   In Section 7, we consider the case of real functions ($X=\R$);
 we show that in this case our conditions easily reduce to those found by Lebedev in \cite{Leb}.

\section{Preliminaries}

By $\lambda$ we will denote the Lebesgue measure on~$\R$ and all integrals are Lebesgue integrals.
Throughout the whole article, $X$ will always be a (real) Banach space.
By $\Hau^1$ we will denote the $1$-dimensional Hausdorff measure. The symbol $\sim$ is used for the strong
 equivalence (i.e., $f \sim g$ means $f/g \to 1$).
\par
A mapping is {\em $L$-Lipschitz} provided it is Lipschitz with some constant $L$
(not necessarily the minimal one).
If $M \subset A \subset \R$ and
$f: A \to X$ are given, then we define the {\em variation of $f$ on $M$} as
\[V(f,M) : = \sup\{\sum^n_{i=1}\|f(x_{i-1})-f(x_i)\|\},\]
where the supremum is taken over all $(x_i)^n_{i=0}\subset M$ such that
$x_0 < x_1 <\dots<x_n$. (We set $V(f,M):=0$, if $M$ is empty or a singleton.)
We say that $f: [a,b] \to X$ is {\em BV}
(or {\em has bounded variation}), provided $V(f,[a,b])<\infty$.

For basic well-known properties of variation, see, e.g., \cite{Fe} and \cite{Chl}. In particular, we will need the additivity
 of variation (see \cite[(P3) on p. 263]{Chl}):
 \begin{equation}\label{adva}
V(f,M) =  V(f, M \cap (-\infty,t]) +   V(f, M \cap [t,\infty)),\ \ \ \text{whenever}\ \ \ t \in M.
 \end{equation}
If $f: [a,b] \to X$ is BV,  then   we define
$v_f(x):=V(f,[a,x]),\ x\in[a,b]$. 
If $f$ is also continuous, then $v_f$ is  continuous  as well (\cite{Fe}, \cite{Chl}). Moreover, clearly $v_f$ is (strictly) increasing,
 if and only if $f$ is not constant on any subinterval of $[a,b]$.
We say that $f: [a,b] \to X$  {\em is parametrized by the arc-length}, if 
 $V(f,[u,v]) = v-u$ for every
$a\leq u <v\leq b$. Obviously, each such $f$ is $1$-Lipschitz (\cite[p. 267]{Chl}).
\begin{definition}\label{arclength}
\begin{itemize}
\item[(a)]
Let  $f: [a,b] \to X$ be a continuous mapping.
We say that $f^*: [c,d] \to X$ is a parametrization of $f$ if
there exists an increasing homeomorphism $h: [c,d] \to [a,b]$ such that $f^* = f \circ h$.
 If $f^*$ is moreover parametrized by the arc-length,
 we say
that $f^*$ is an {\it arc-length parametrization of $f$}.
\item[(b)]
If $f: [a,b] \to X$ is  nonconstant, continuous and BV, then
there exists (see \cite[\S2.5.16]{Fe} or \cite[Theorem 3.1]{Chl})  a unique $F: [0,\ell] \to X$ (where $\ell:=v_f(b)$)
such that $f = F \circ v_f$.
We will denote this associated mapping $F$ by $\mcA_f$.
\end{itemize}
\end{definition}
It is easy to see that $\mcA_f$ is always parametrized by the arc-length, and thus it is $1$-Lipschitz (see \cite{Chl}). 
We will use several times the following easy lemma.
\begin{lemma}\label{alp}
Let  $f:[a,b] \to X$ be continuous. Then the following hold.
\begin{enumerate}
\item
 The function $f$ has an arc-length parametrization if and only if
$f$ is BV and $f$ is not constant on each $[c,d]\subset [a,b]$.
In this case, $\mcA_f$ is  an arc-length parametrization of $f$,  $\mcA_f= f \circ (v_f)^{-1}$,
 and a general  arc-length parametrization of $f$ is of the form $F^s (x) = \mcA_f(x-s),\ x \in [s,s+ \ell],$
  where $\ell:=v_f(b)$ and  $s \in \R$.
\item  
If $f$ is  BV on $[a,b]$, and is not constant on each
subinterval of an interval $[\alpha, \beta] \subset [a,b]$, 
then $\mcA_f|_{[v_f(\alpha),v_f(\beta)]} = f \circ (v_f|_{[\alpha,\beta]})^{-1}$ is an arc-length
parametrization of $f|_{[\alpha,\beta]}$.
\end{enumerate}
\end{lemma} 
Let $f:[a,b]\to X$. The derivative $f'$ is defined in the usual way; at the endpoints we take
the corresponding unilateral derivatives.
We say that {\em $f: [a,b] \to X$ is $C^{n}$} ($n\in\N$) provided the $n$-th derivative $f^{(n)}$ exists
and is continuous on $[a,b]$.
We will say that {\em $f: [a,b] \to X$ is $D^{2,\infty}$} if $f''$ exists and is bounded on $[a,b]$.
Clearly, if $f$ is $C^2$, then $f$ is $D^{2,\infty}$. Further, if $f$ is $D^{2,\infty}$, then $f'$ is clearly Lipschitz. 

We will need also the following almost obvious lemma.
\begin{lemma} \label{eqnorm}
Let $X$ be a Banach space and $f:[0,1] \to X$. Then the validity of the statement that $f$ is $C^2$ (resp. $D^{2,\infty}$) does not depend on a choice of an equivalent norm on $X$.
\end{lemma}

It is well known (see e.g. \cite[p.2]{VZ}, or use \cite[Theorem~7]{K} together with~\cite[Theorem~2.10.13]{Fe})
that if $f:[a,b]\to X$ is Lipschitz, and $f'(x)$ exists for almost all $x\in[a,b]$,
then 
\begin{equation}\label{varder} 
V(f,[a,b])=\int^b_a\|f'(x)\|\,dx.
\end{equation}
\par

For a proof of the following well-known version of Sard's
Theorem, see e.g.\ \cite[Theorem~7]{K}.

\begin{lemma}\label{Sard} 
Let $f:[0,1]\to X$ be arbitrary.
Let $C: =  \{x\in[0,1]:\ f'(x)=0\}$. Then
${\mathcal{H}}^1(f(C))=0.$
\end{lemma}

We will need also the following known lemmas.

\begin{lemma}\label{arclenlem}
If $X$ is a Banach space, $f:[a,b]\to X$ is continuous, BV, not constant on any interval, and such that $F=f\circ (v_f)^{-1}$ is
$C^1$, then $\|F'(x)\|=1$ for each $x\in[0,\ell]$, where $\ell=v_f(b)$.
\end{lemma}

\begin{proof}
Since $F$ is an arc-length parametrization of $f$, and thus $1$-Lipschitz,
we obtain $\|F'(s)\| \leq 1,\ s \in [0,\ell]$, and $\ell =V(F,[0,\ell])=\int^\ell_0 \|F'(s)\|\,ds$ by~\eqref{varder}.
Since $F$ is $C^1$, we obtain $\|F'(s)\|=1$ for each $s\in[0,\ell]$.
\end{proof}

\begin{lemma}\label{H1L}
Let $f:[a,b] \to X$ be continuous. Let $G\subset (a,b)$ be an open set,
$H:= [a,b] \setminus G$ and  $(a_t,b_t),\ t\in T$, be all $($pairwise
different$)$ components of $G$. Then:
\begin{enumerate}
\item  If ${\mathcal H}^1(f(H))=0$, then $V(f,[a,b])= \sum_{t\in T}
 V(f,[a_t,b_t])$.
\item If $V(f,[a,b])= \sum_{t\in T}
 V(f,[a_t,b_t]) < \infty$, then   ${\mathcal H}^1(f(H))=0$.
\item If    ${\mathcal H}^1(f(H))=0$ and $f$ is $L$-Lipschitz  on
each $[a_t,b_t]$, then $f$ is $L$-Lipschitz on $[a,b]$.
\item If $f$ is BV and  ${\mathcal H}^1(f(H))=0$, then 
 $\lambda(v_f(H)) = 0$.
 \end{enumerate}
\end{lemma}

\begin{proof}
Part~(i) is an easy consequence of the vector form (\cite[Theorem 2.10.13]{Fe})
of Banach indicatrix theorem; see ~\cite[Lemma~2.7]{DZ}. For  part~(ii), let $\e>0$ and choose a finite $S\subset T$
 such that
$\sum_{t\in S} V(f,[a_t,b_t])>V(f,[a,b])-\e$.
Then by~\cite[Corollary~2.10.12]{Fe} we see that
\[ \Hau^1(f(H)) \leq \sum^N_{i=1} \Hau^1(f([c_i,d_i]))\leq \sum^N_{i=1} V(f,[c_i,d_i])\leq \e,\]
where $[c_i,d_i]$, $i=1,\dots,N$, are components of
$[a,b] \setminus \bigcup_{t\in S} (a_t,b_t)$.
Therefore $\Hau^1(f(H))=0$.
\par
For part~(iii), consider arbitrary  $c,d \in H$, $c<d$. 
Applying~(i) to $f$ on $[c,d]$, we obtain
\begin{equation*}
\begin{split} 
|f(d)-f(c)| &\leq V(f,[c,d]) = \sum \{V(f,[a_t,b_t]): t\in T,
[a_t,b_t] \subset [c,d]\} \\
&\leq \sum \{L |b_t-a_t|: t\in T,[a_t,b_t] \subset [c,d]\} \leq L(d-c),
\end{split}
\end{equation*}
where we used the fact that the variation of an $L$-Lipschitz function
on an interval is also $L$-Lipschitz.
Since $f$ is $L$-Lipschitz on $H$ and on each $[a_t,b_t]$, it is clearly
$L$-Lipschitz on $[a,b]$. Part (iv) easily follows from (i).
\end{proof}

For the following well-known fact see \cite[Theorem 2.3, p. 35]{B}.
\begin{lemma}\label{bru}
Let $C\subset [0,1]$ be a closed nowhere dense set. Then there exists a real function $\vf$ on $[0,1]$ which has a bounded derivative on $[0,1]$, and $C$ is the set of points of discontinuity of $\vf'$.
\end{lemma}

\section{Basic special notions and lemmas}

The following easy inequality is well known (see e.g.\ \cite[Lemma~5.1]{MS}):

\begin{equation}\label{triangle}
\text{if }u,v\in X\setminus\{0\},\text{ then }
 \left\|\frac{u}{\|u\|}-\frac{v}{\|v\|}\right\|
\leq \frac{2}{\|u\|}\, \|u-v\|.
\end{equation}

\begin{lemma}
\label{lipkonst} 
Let $I$ be a closed interval and $f:I\to X$ be such  that $\|f'(x)\|\geq i>0$ and
 $\|f''(x)\| \leq M$ for each $x \in I$. Let $F$ be an arc-length parametrization of $f$,
  and suppose that  $F''(z)$ exists.
 Then
 $\|F''(z)\|\leq\frac{2\,M}{i^2}$.  
\end{lemma}

\begin{proof}
 By Lemma \ref{alp}(i), we can suppose that
$F=f\circ v_f^{-1}$.
Let $x,y\in I$. Since $f'$ is $M$-Lipschitz, by~\eqref{triangle} we have
\[ \left\|\frac{f'(x)}{\|f'(x)\|}-\frac{f'(y)}{\|f'(y)\|}\right\|
\leq \frac{2}{\|f'(x)\|}\|f'(x)-f'(y)\|
\leq \frac{2\,M}{i}|x-y|.\]
Consequently, $g:=\frac{f'}{\|f'\|}$ is $\frac{2M}{i}$-Lipschitz.
Since~\eqref{varder} implies $(v_f)' = \|f'\|$, we have  $(v_f^{-1})' = \frac{1}{\|f' \circ v_f^{-1}\|}$.
Thus, $v_f^{-1}$ is $\frac{1}{i}$-Lipschitz; so
  $F' = (f \circ v_f^{-1})' = g \circ v_f^{-1}$ is $\frac{2M}{i^2}$-Lipschitz on $v_f(I)$, and
   the conclusion 
follows.
\end{proof}

We say that a Banach space {\em $X$ has a Fr\'echet smooth norm}, provided 
the norm $\|\cdot\|$ of $X$ is Fr\'echet differentiable at all $x\in X\setminus\{0\}$
(or equivalently $\|\cdot\|$ is $C^1$ on $X\setminus\{0\}$).

\begin{lemma}
\label{vlastF}
Let $X$ be a Banach space which admits an equivalent Fr\'echet smooth norm. 
 Let $I$ be a closed interval and $f:I\to X$ be such  that $f'(x) \neq 0$ for each $x \in I$.
 If $f$ is $D^{2,\infty}$ (resp. $C^2$), then 
   $F:=f\circ v_f^{-1}$ is $D^{2,\infty}$ (resp. $C^2$).
\end{lemma}
\begin{proof} We can assume that $X$ has a Fr\'echet smooth norm by Lemma \ref{eqnorm}.
Let $f$ be $D^{2,\infty}$. Since clearly $i:= \min\{\|f'(x)\|:\ x \in I\} >0$, by Lemma \ref{lipkonst}
 it is sufficient to observe that $F''(y)$ exists for each $y \in v_f(I)$. In the proof of Lemma \ref{lipkonst}
  we obtained the equality  $(v_f^{-1})' = \frac{1}{\|f' \circ v_f^{-1}\|}$, which implies that $(v_f^{-1})'$
   is differentiable. So, $v_f^{-1}$ is twice differentiable, which implies that $F:=f\circ v_f^{-1}$ is 
   twice differentiable. 
   
   If $f$ is $C^2$, then the above equality implies that $(v_f^{-1})'$ is $C^1$. So, $v_f^{-1}$ is $C^2$,
    which implies that $F:=f\circ v_f^{-1}$ is $C^2$. 
\end{proof}

Now we define a special type of an ``1/2-variation'', which is crutial in our solution of $C^2$ and $D^{2,\infty}$ parametrization problems.

\begin{definition}\label{dvojve} 
Let $f:[0,1]\to X$ be continuous and BV. Let
$\emptyset\neq G \subset (0,1)$ be an open set and $0  < \delta < \infty$.
Suppose that $f$ is not constant on each interval contained in $G$ and 
$F:=\mcA_f$
(see Definition~\ref{arclength}) is twice differentiable on $v_f(G)$. Then we define
\[ W^{\delta}(f,G) = \sup \big\{ \sum_{k=1}^{n}\ \sqrt{V(f,I_k)}\, \big\},\]
where the supremum is taken over all non-overlapping systems
$I_1,\dots,I_n$ of compact intervals with $\inte(I_k) \subset G$ such
that $S_{I_k}\cdot V(f,I_k) \geq \delta$ whenever  $I_k \subset G$
 (where $S_{I_k}=\sup_{x\in v_f(\inte(I_k))} \|F''(x)\|\leq \infty$).
\end{definition}

\begin{remark}\label{Wdeltainv}
Let $f$, $G$ and $\delta>0$ be as in Definition \ref{dvojve}, and let 
$\omega:[0,1]\to[0,1]$ be an increasing homeomorphism. Then
\[W^\delta(f,G)=W^\delta(f\circ\omega,\omega^{-1}(G)).\] 
This equality easily follows from the definition, if
 we observe that $\mcA_f = \mcA_{f \circ \omega}$ and $v_{f \circ \omega} = v_f \circ \omega$.
\end{remark}

\begin{remark}\label{komponenty}
Let $f$, $G$ and $\delta>0$ be as in Definition \ref{dvojve}, 
and let $\mcI$ be the family   of all components of $G$. Then, using only Definition \ref{dvojve}, we clearly obtain
$$\sum_{I\in\mcI} \sqrt{V(f,I)} \leq W^\delta(f,G). $$.
\end{remark}

Our first basic lemma follows.

\begin{lemma}\label{Wfin2}
Suppose that $f: [0,1] \to X$ is continuous and BV, and  $f''$ exists and is bounded on $(0,1)$.  
Let $\emptyset \neq G \subset (0,1)$ be an
open set such that $f'(e)=0$, whenever $e\in (0,1)$ is an
endpoint of any component of $G$. Let $f$ be nonconstant on each interval contained in $G$ 
and let $F:=\mcA_f$ $($see Definition~\ref{arclength}$)$ be twice differentiable on
 $v_f(G)$. Then, $W^\delta(f,G) < \infty$ for each $\delta>0$.
\end{lemma}

\begin{proof}
Denote $M := \sup_{x \in (0,1)} \|f''(x)\|$. Choose $\delta>0$  and
 consider a system $I_1,\dots,I_n$ 
of non-overlapping compact intervals with $\inte(I_k) \subset G$ such
that $S_{I_k}\cdot V(f,I_k) \geq \delta>0$ whenever  $I_k \subset G$
(where $S_{I_k}=\sup_{x\in v_f(\inte(I_k))} \|F''(x)\|$).

Further, consider an arbitrary $I_k =: I $ such that $I \subset (0,1)$. Denote  $i:=\min\{\|f'(x)\|:x\in I\}$
 and choose $\tilde x \in I$ such that $\|f'(\tilde x)\| = i$. The Mean Value Theorem implies that, for each
  $x \in I$,
  \[ \|f'(x)\| - i \leq \|f'(x) - f'(\tilde x)\| \leq M \lambda(I).\]
  So, $\|f'(x)\| \leq i + M \lambda(I)$, and thus \eqref{varder} implies
  \begin{equation}\label{odhvar}
  V(f,I) = \int_{I} \|f'(x)\|\ dx \leq i \lambda(I) + M (\lambda(I))^2.
  \end{equation}
  Thus, if $i=0$, we obtain
  \begin{equation}\label{ject}
  \sqrt{V(f,I)} \leq \sqrt{M}\, \lambda(I).
  \end{equation}
  
  If $i>0$, then the assumptions of the lemma imply that $I \subset G$, and therefore $\sup_{x\in v_f(\inte(I))} \|F''(x)\| \, V(f,I) \geq \delta$.
By Lemma \ref{lipkonst} and Lemma \ref{alp}(ii), we obtain
\begin{equation}\label{odhi}
\frac{2M}{i^2}\, V(f,I) \geq \delta,\ \ \ \text{and so}\ \ \ i \leq \sqrt{\frac{2M}{\delta} V(f,I)}. 
\end{equation} 
 Using \eqref{odhvar}, \eqref{odhi} and the A-G inequality, we obtain subsequently
 \[ V(f,I) \leq  \sqrt{\frac{2M}{\delta} V(f,I)}\, \lambda(I) + M (\lambda(I))^2 = 
  \sqrt{V(f,I) \cdot \frac{2M}{\delta}(\lambda(I))^2} +   M (\lambda(I))^2,\]
\[ V(f,I) \leq \frac{V(f,I)}{2} + \frac{M}{\delta}(\lambda(I))^2 +  M (\lambda(I))^2,\]
\begin{equation}\label{dvact}
\sqrt{V(f,I)} \leq \sqrt{2\left(\frac{M}{\delta} + M\right)} \, \lambda(I).
\end{equation}
Since there are at most two intervals $I \in \{I_1,\dots,I_n\}$ which are not contained in $(0,1)$, and 
 $\lambda(I_1)+\cdots + \lambda(I_n) \leq 1$, we obtain
 by  \eqref{ject} and \eqref{dvact}  
\[ \sum_{k=1}^n \sqrt{V(f,I_k)} \leq \sqrt{M} + \sqrt{2\left(\frac{M}{\delta} + M\right)} + 2 \sqrt{V(f,[0,1])}.\]
Consequently,
$W^{\delta}(f,G) < \infty$.
\end{proof}

We will work many times with the following natural notion of a generalized partition of an open set.

\begin{definition}\label{nekdel}
\begin{enumerate}
\item We say that
 $\mcI\subset\Z$ is an $\Z$-interval, if
$\mcI=(l,m) \cap \Z $, where $l,m\in\Z^*=\Z\cup\{-\infty,\infty\}$.
\item We will say that a family $\p$ of compact intervals is a
{\it generalized partition of a bounded interval $(a,b)$}, if there exists
 a system  $(x_i)_{i\in\mcI}$ such that $\mcI$ is an $\Z$-interval,
 the function $i \mapsto x_i,\ i \in \mcI$, is strictly increasing,
 $\inf_{i\in\mcI} x_i=a$,  $\sup_{i\in\mcI} x_i=b$ and
$\p = \{ [x_k,x_{k+1}]:\ k, k+1 \in \mcI\}$.
\item  We will say that a family $\p$ of compact intervals is a
{\it generalized partition of a bounded open set $\emptyset \neq G \subset \R$}, if
$\bigcup \{\inte(I): I \in \p\} \subset G$, and for each component
$(a,b)$ of $G$, the family $\{I \in \p:\ \inte(I) \subset (a,b)\}$
is a generalized partition of $(a,b)$.
\end{enumerate}
\end{definition}

In fact, we will work only with special generalized partitions. In the $D^{2,\infty}$ problem, we need the notion
 of an $[f,\delta,K]$-partition. In the $C^{2,\infty}$ problem, it is sufficient to work only with the more special (and simpler) notion of an $(f,\delta,K)$-partition.

\begin{definition}\label{specdel}
Suppose that   $f: [0,1] \to X$  is continuous and BV. 
Further suppose that $f$ is not constant on any interval which is a subset of an open set $\emptyset \neq G \subset (0,1)$. Set 
$F:=\mcA_f$ (see Definition~\ref{arclength}) and suppose that $F''$ exists on
$v_f(G)$.  For each compact interval $L$ with $\inte(L) \subset G$, set $S_L:=  
 \sup_{x\in v_f(\interior(L))}\|F''(x)\|$.
 Let $\p$ be a generalized
partition of $G$, $0\leq \delta \leq K \leq \infty$, and $\delta\in\R$. Then we say
that  $\p$ is an {\it  $[f,\delta,K]$-partition of $G$}, if:
\begin{enumerate}
\item[(a)]
$S_I<\infty$\ \ and \ \  $S_I\cdot V(f,I) \leq K$\ \  for each\ \  $I \in \p$.
\item[(b)]
For each $I \in \p$ with $I \subset G$, there exists $J \in \p$ such that $I \cap J \neq \emptyset$
 and $S_{I \cup J} V(f,I\cup J) \geq \delta$.
\end{enumerate}
We say
that  $\p$ is an {\it  $(f,\delta,K)$-partition of $G$}, if (a) holds, and
\begin{enumerate}
\item[(b*)]
  $S_I\cdot V(f,I) \geq \delta$
for each $I \in \p$ with $I \subset G$.
\end{enumerate}
\end{definition}

\begin{remark}\label{kulhr}
Clearly each $(f,\delta,K)$-partition of $G$ is an $[f,\delta,K]$-partition of $G$.
 The notions of an $[f,0,K]$-partition and of an $(f,0,K)$-partition coincide. 
\end{remark}

\begin{lemma}\label{konstrpartlem} 
Let $\emptyset\neq G\subset(0,1)$ be open, 
and $f:[0,1]\to X$ be a continuous BV function.
Let $f$ be not constant on any subinterval of $G$, and let
$F:=\mcA_f$ $($see Definition~\ref{arclength}$)$ have locally bounded second derivative (resp. be 
$C^2$) on $v_f(G)$. Then, for each  $ 0 <\delta < \infty$,
 there exists a generalized partition $\mcP$ of $G$,
which is an  $[f,\delta,\delta]$-partition (resp. $(f,\delta,\delta)$-partition) of $G$.
\end{lemma}

\begin{proof} 
Without any loss of generality, we can assume that $G=(a,b)\subset (0,1)$.
Let $x_0:=\frac{a+b}{2}$.
We will construct points $x_i$ ($i\in\Z$) with
\[
a \leq \dots \leq x_{-2} \leq x_{-1} < x_0 < x_1 \leq x_2 \leq x_3 \leq \dots \leq b.
\]
Suppose that  $x_{n-1}$ ($n \in \N$) is defined. If $x_{n-1}=b$, then put $x_n = b$.
 If  $x_{n-1}<b$, then 
 set  $x_n := \inf M_n$, where 
\[M_n:=\{t\in (x_{n-1},b]:\ V(f,[x_{n-1},t])
\cdot   \sup_{x\in v_f(x_{n-1},t)}\|F''(x)\|   \geq\delta\}\ \cup \{b\}.\]
 Since  $F''$ is locally bounded on $(a,b)$, we easily see that $x_n > x_{n-1}$. Further, it is easy to see
 that
\begin{equation}\label{lt1}
 V(f,[x_{n-1},x_n])
\cdot   \sup_{x\in v_f(x_{n-1},x_n)}\|F''(x)\|      \leq\delta.
\end{equation}
 (Otherwise , using the definition of the supremum, we easily see that there exists a $t \in M_n \cap (x_{n-1},x_n)$.)
\par
We define the points $x_n$ ($n < 0$) quite symmetrically and set 
 $\p := \{[x_i,x_{i+1}]:\ i \in \Z, x_i < x_{i+1}\}$. 
Since  $F''$ is locally bounded on $(a,b)$, it is easy to show that $\inf_{n\in\mcJ} x_n=a$ and $\sup_{n\in\mcJ} x_n=b$. So, $\p$ is a generalized partition
 of $(a,b)$. 
To prove that $\p$ is an $[f,\delta,\delta]$-partition of $(a,b)$, first observe that \eqref{lt1} holds 
 whenever $[x_{n-1},x_{n}] \in \p$ (also if $n-1<0$), and so the condition (a) of Definition \ref{specdel}
  holds. To prove also (b), consider an interval $I \in \p$ such that $I \subset (a,b)$, i.e.
  $x_{n-1} \neq a$ and $x_{n}\neq b$. 
If $n-1\geq0$,
then clearly $x_{n+1} \in M_n$. So, putting $J:= [x_n,x_{n+1}]$, we obtain $S_{I \cup J} V(f,I\cup J) \geq \delta$.
 Similarly, if $n-1 < 0$, we can set $J:= [x_{n-2},x_{n-1}]$.
\par
In the case that $F$ is $C^2$ on $v_f(G)$, then $\p$ is even an $(f,\delta,\delta)$-partition of $(a,b)$.
 Indeed, let  $I = [x_{n-1},x_{n}]$, $n-1\geq 0$, and $x_n \neq b$. Using continuity of $F''$ in $v_f(x_n)$, we easily
  obtain $S_{I}\cdot  V(f,I) = \delta$. Using a symmetrical argument, we obtain that the condition
   (b*) of Definition \ref{specdel} holds. 

\end{proof}

\begin{lemma}\label{nero}
Let $a_i$ $(i \in I)$, $b_j$, $c_j$ $(j\in J)$ be non-negative numbers,
$I$ countable, and $J$ finite. 
Then
\begin{equation}\label{dvenero}
\sqrt{ \sum_{i\in I} a_i} \leq \sum_{i \in I} \sqrt{a_i}\ \ 
\mbox{and}\ \ 
\sum_{j\in J} \sqrt{b_j c_j} \leq \sqrt{ \sum_{j\in J} b_j \cdot\sum_{j\in J} c_j}.
\end{equation}
\end{lemma}

\begin{proof}
The first inequality is clear. The second is an immediate
consequence of the Cauchy-Schwartz inequality.
\end{proof}

\begin{lemma}\label{nove}
Let $f:[0,1]\to X$ be continuous and BV. Let
$\emptyset\neq G \subset (0,1)$ be an open set and $0  < \delta < \infty$.
Let $f$ be  nonconstant on each interval contained in $G$, let
$F:=\mcA_f$
(see Definition~\ref{arclength}) have locally bounded second derivative on $v_f(G)$, and let
 $W^{\delta}(f,G)<\infty$.
Then the following hold.
\begin{enumerate}
\item[(i)]
If $\p$ is an $[f, \delta, \infty]$-partition
of $G$, then   $ \sum_{I \in \p} \sqrt{V(f,I)} < \infty$. 
 \item[(ii)] 
 $\int_{v_f(G)}\sqrt{\|F''\|} < \infty$.
\end{enumerate}
\end{lemma}

\begin{proof}
Let $\p$ be an $[f, \delta, \infty]$-partition
of $G$. For each compact interval $L$ with $\inte(L) \subset G$, set $S_L:=  
 \sup_{x\in v_f(\interior(L))}\|F''(x)\|$.
Set $\p'=\{I\in\p:I\subset G\}$ and $\p''=\p\setminus\p'$.
By Definition~\ref{dvojve}, we have 
\begin{equation}\label{zdef}
\sum_{I\in\p''} \sqrt{V(f,I)}\leq W^\delta(f,G) < \infty
\end{equation}
By Definition~\ref{specdel}, for each $I\in\p'$ 
we can choose an interval $J =:n(I)$  such that  
\[ V(f,I\cup n(I))\cdot S_{I\cup n(I)}\geq\delta.\]
It is easy to see that $\{I\cup n(I):I\in\p'\}$ can be written as $\bigcup^3_{i=1}\p_i$,
where $\p_i$ is a family of non-overlapping compact intervals for each $i\in\{1,2,3\}$.
Thus,
\[ \sum_{I\in\p'} \sqrt{V(f,I)}\leq
\sum_{i=1}^3 \sum_{I\in\p_i} \sqrt{V(f,I\cup n(I))} \leq
3 W^\delta(f,G) < \infty.\]
So, 
\eqref{zdef} implies $ \sum_{I \in \p} \sqrt{V(f,I)} < \infty$.
Thus, we have proved (i).

To prove (ii), choose by Lemma~\ref{konstrpartlem} an $[f,\delta,\delta]$-partition $\p$ of $G$. 
 Observe that $\|F''\|$ is Lebesgue measurable. Thus
\begin{equation*}
\begin{split}
&\int_{v_f(G)}\sqrt{\|F''\|}\leq\sum_{I\in\p}\lambda(v_f(I))\cdot \sup_{x \in v_f(\inte(I))}\sqrt{\|F''(x)\|}\\
&=\sum_{I\in\p}\sqrt{V(f,I)}\cdot\sup_{x \in v_f(\inte(I))}\sqrt{\|F''\|\cdot V(f,I)}
\leq\sqrt{\delta}\sum_{I\in\p}\sqrt{V(f,I)}<\infty,
\end{split}
\end{equation*}
where the last inequality follows from (i).
\end{proof}

\begin{lemma}\label{vt2}
Let  $f: [0,1] \to X$ be continuous and BV. Let $\emptyset \neq G \subset (0,1)$ be an open set
 such that $f$ is nonconstant on any interval contained in $G$ and let
$F:=\mcA_f$ (see Definition~\ref{arclength}) have locally bounded second derivative  on $v_f(G)$.
Suppose that  $\s$ is a family of pairwise non-overlapping compact intervals
such that $\inte(J) \subset G$ for each $J \in \s$ and
\begin{equation}\label{tri22}
\sum_{J\in \s} V(f,J)= V(f,[0,1]),\  \sum_{J\in \s} \sqrt{V(f,J)} <
\infty,\
\sum_{J\in \s}  V(f,J) \sqrt{S_J} < \infty,
\end{equation} 
where $S_J=\sup_{y\in v_f(\inte(J))}\|F''(y)\|$.
Then $W^{\delta}(f,G) < \infty$ for each $\delta > 0$.
\end{lemma}

\begin{proof}
Let $\delta >0$ and consider a finite system $\Ka$
of non-overlapping compact intervals such that $\inte(I) \subset G$ for each $I \in \Ka$, and $S_{I}\cdot V(f,I) \geq \delta$ 
whenever  $I \in \Ka$ and $I \subset G$. For each $J \in \s$, let  
$\Ka_J : = \{ I \in \Ka: I \subset \inte(J)\}$. 
Set $\Ka_1 := \bigcup\{\Ka_J: J \in \s\}$ and $\Ka_2 : = \Ka \setminus \Ka_1$.    
For each $J \in\s$, we obtain 
\[\sqrt{\delta} \sum\{\sqrt{V(f,I)}: I\in \Ka_J\} \leq
\sum \{\sqrt{S_I}\cdot V(f,I): I\in \Ka_J\} \leq
\sqrt{S_J}\cdot V(f,J).\]
Therefore
\begin{equation}\label{K1c}
\sum\{\sqrt{V(f,I)}: I \in \Ka_1\} \leq (1/\sqrt{\delta}) \sum\{\sqrt{S_J}\cdot V(f,J):
J \in \s\}.
\end{equation}
\par
For each $I \in \Ka_2$, denote by $\s_I$ the set of all $J\in \s$, such
that $J \cap \inte(I) \neq \emptyset$.
If  $I = [a,b] \in \Ka_2$, put $a^* := \min J_a$, if there exists $J_a \in \s$
with $a \in \inte(J_a)$ and $a^* := a$, if such $J_a$ does not exist. Similarly,
put $b^* := \max J_b$, if there exists $J_b \in \s$
with $b \in \inte(J_b)$ and $b^* := b$, if such $J_b$ does not exist.
    The equality of~\eqref{tri22} easily implies that 
$ V(f, [a^*,b^*]) = \sum \{V(f,J): J\in \s_I\}$. (It can be shown either directly, or using first  Lemma~\ref{H1L}(ii)
 and then Lemma~\ref{H1L}(i).)
Thus
the first inequality of ~\eqref{dvenero} implies   
$\sqrt{V(f, [a^*,b^*])} \leq \sum \{\sqrt{V(f,J)}: J \in \s_I\}$. 
Observing that, for each $J \in \s$, the set 
$\{I\in \Ka_2:J \in \s_I\}$ contains at most two intervals, we obtain
\begin{equation}\label{K2c}
\sum\{\sqrt{V(f,I)}: I \in \Ka_2\} \leq 2 \sum\{\sqrt{V(f,J)}:J \in \s\}.
\end{equation}
\par
Now~\eqref{tri22},~\eqref{K1c} and~\eqref{K2c} imply $W^{\delta}(f,G) < \infty$.
\end{proof}

\begin{lemma}\label{homint}
Let $ (I_{\alpha})_{\alpha\in A}$ be a system of
pairwise non-overlapping compact subintervals of an interval~$[0,d]$. Let
$\sum_{\alpha \in A} \mu_{\alpha} < \infty$, 
where $\mu_{\alpha}>0$, $\alpha \in A$. Then there exists an interval $[0,d']$ and 
an increasing homeomorphism $\Psi: [0,d'] \to [0,d]$ such that $\lambda(\Psi^{-1}(I_{\alpha})) = \mu_{\alpha}$
and $\Psi^{-1}$ is absolutely continuous.
\end{lemma}

\begin{proof}
We can define $\Psi := \omega^{-1}$, where
$\omega(x) = \int_0^x \vf$ ($x \in [0,d']$), and
$\vf(t) = \mu_{\alpha}/\lambda(I_{\alpha})$ for $t \in \interior(I_{\alpha})$ and
$\vf(t) = 1$ for $t \in [0,d'] \setminus \bigcup\{\interior(I_{\alpha}): \alpha\in A\}$.
\end{proof}

\begin{lemma}\label{l1}
Let $0<d<1$, $\xi > 0$ and $c_l, c_r \geq  0$ with $0 < \max (c_l,c_r) \leq
10^{-9}\cdot\xi d$ be given. Then, for each  $I = [u,v]$
with $\lambda(I) = d$, there exists a $C^2$ function $\omega$ on $I$ such that
$\omega(u) = \omega(v) = 0$, $\omega'(x) = c_l$ for $x \in [u,u+d/3]$, $\omega'(x) =
c_r$ for $x \in [u+2d/3,v]$, and $\max (|\omega'(x)|, |\omega''(x)|) \leq \xi$ for $x\in I$.
\end{lemma}

\begin{proof}
We can suppose that $I= [0,d]$. First assume that $c_l \leq c_r$.
Denote $p_1 := 2 (9 c_r -c_l) \xi^{-1}$, $p_2 := 20\cdot c_r\cdot
\xi^{-1}$; clearly
\begin{equation}\label{pjpd}
0 < p_1 \leq p_2 < 10^{-7} d.
\end{equation}
For every  $0 \leq s \leq d/3 - p_1 - p_2$,  let $x_0 := 0$, $x_1 :=
d/3$, $x_2 := x_1 + p_1/2$, $x_3 := x_2 +  p_1/2$, $x_4 := x_3 + s$,
$x_5 := x_4 + p_2/2$, $x_6 := x_5 + p_2/2$, $x_7 := d$.
Clearly $x_i \leq x_{i+1}$ ($0\leq i \leq 6$). Now consider the 
function $\mu = \mu_s$ which is linear on each $[x_i,x_{i+1}]$ (if $x_i
< x_{i+1}$) and
\[\mu(x_0) = \mu(x_1) = \mu(x_3) = \mu(x_4) = \mu(x_6) = \mu(x_7)=0,\
\mu(x_2) = -\xi,\  \mu(x_5) = \xi.\]
Clearly  $|\mu(x)| \leq \xi$, $x \in [0,d]$.
Let  $\nu(x) = \nu_s(x)= c_l + \int_0^x \mu$, $x\in I$. It is easy to
see that $\nu(x) = c_l$ for $x \in [x_0,x_1]$, $\nu(x) = -9c_r$ for $x \in [x_3,x_4]$, 
$\nu(x)=c_r$ for $x \in [x_6,x_7]$, and $|\nu(x)| \leq 9 c_r \leq \xi$ for $x \in I$.
These properties of $\nu$ and~\eqref{pjpd} easily imply that
$g(s) := \int_0^d \nu_s >0$ for $s=0$ and $g(s) < 0$ for
$s = d/3 -p_1 - p_2$. Since $g(s)$ is clearly continuous, we can choose
$s_0 \in (0, d/3 -p_1 - p_2)$ with $g(s_0)=0$. 
Now it is easy to see that $\omega(x) :=\int_0^x \nu_{s_0}$, $x \in I$, 
has all desired properties.
\par
If $c_l > c_r $, we apply the just-proven assertion to $I^* :=
[-d,0]$, $c_l^* := c_r$, $c_r^* := c_l$ and obtain a function $\omega^*$
on~$I^*$. Now it is sufficient to put $\omega(x) := -\omega^*(-x)$, $x \in [0,d]$.
\end{proof}

\begin{lemma}\label{l2}
Let  $V>0 $, $ \eta > 0$ and $c_l, c_r \geq 0$ such that
$d:= \sqrt{V/\eta} < 1$ and $\max(c_l,c_r) \leq 10^{-10} \cdot d^2\eta$ be given. 
Then, for every interval $J$ of the length $V$ and
every interval $I = [u,v]$ of the length $d$, there exists an increasing $C^2$ 
homeomorphism $\vf : I \to J$ such that $\vf'(u)=c_l$, $\vf'(v) = c_r$,
$\vf''(u)=\vf''(v)=0$,   $0 < \vf'(x) \leq  19\cdot \sqrt {\eta V}$ for $x \in \inte(I)$,
and $|\vf''(x)| \leq  19\cdot \eta$ for $x \in I$.
\end{lemma}

\begin{proof}
We can suppose that  $J = [0, V]$ and $I = [0,d]$.  
Consider the function $\pi$ on $I$ which is linear on each $[id/6, (i+1)d/6]$
($i=0,\dots,5$), $\pi(d/6) = \eta$, $\pi (5d/6)= -\eta$ and
$\pi(id/6) =0$ for $i=0,2,3,4,6$. Define  $\rho(x) := \int_0^x \pi$,
$x \in I$, and $\tau(x) := \int_0^x \rho$, $x \in I$. 
Clearly $\rho(0)= \rho(d) = 0$,
$0 < \rho(x) \leq \eta d$ for $x \in (0,d)$ and
$\rho(x) = d\eta/6$ for $x \in [d/3, 2d/3]$. 
Therefore $ d^2\eta/18 \leq \tau(d) \leq \eta d^2$. 
Thus, setting $\Psi(x) = (\eta d^2/\tau(d)) \tau(x)$ ($x \in I$), we have
$\Psi(d) = \eta d^2 = V$. 
Consequently, $\Psi: I \to J$ is an increasing homeomorphism with $\Psi'>0$ on $\inte(I)$.
Further
\begin{equation}\label{Psi}
|\Psi''(x)| \leq 18 \cdot |\tau''(x)| \leq 18\cdot \eta,\ \ \
\Psi'(x) \leq 18 \cdot \tau'(x) \leq 18\cdot \eta d = 18 \sqrt{\eta V}\ \ \ ( x \in I), \
\end{equation}
\begin{equation}\label{prost}
\Psi'(0) = \Psi'(d)=0,\ \ \text{and}\ \ \    \Psi'(x) \geq \tau'(x) = \rho(x) = \frac{d \eta}{6}, \ \ x \in
\left[\frac{d}{3}, \frac{2d}{3}\right].
\end{equation}
Therefore, in the case $c_l = c_r = 0$, the function $\vf := \Psi$ has the desired properties.
\par
If $\max (c_l,c_r) > 0$, set  $\xi: = d\eta/7$. Then
$\max(c_l,c_r) \leq 10^{-10} \cdot d^2\eta  \leq 10^{-9}\cdot  \xi d $. 
Thus we can use Lemma \ref{l1} and obtain a corresponding function $\omega$ on $[0,d]$.
Define $\vf: = \Psi + \omega$. 
Since $|\omega'(x)| \leq d\eta/7$
if $x \in [d/3, 2d/3]$ and $\omega'(x) \geq 0$ otherwise, using ~\eqref{prost} we obtain 
 that $\vf' >0$ on $\inte(I)$. 
Since $\omega(0) = \omega(d) =0$, we
obtain that $\vf$ is a homeomorphism of $I$ onto $J$.
Since  $\omega'(0)=c_l$, $\omega'(d)=c_r$,  $\max (|\omega'(x)|, |\omega''(x)|) \leq \xi$ for $x\in I$,
 and $\xi = \sqrt{V \eta}/7 = d\eta/7 < \eta$, we obtain by \eqref{Psi} and \eqref{prost} that
 $\vf$ has all desired properties.
\end{proof}

\begin{lemma}\label{hlavkon2}
Suppose that  $f: [0,1] \to X$ is BV continuous, $\emptyset \neq G \subset (0,1)$
is an open set, $f$ is not constant on any subinterval of $G$ and $F:=\mcA_f$ $($see Definition~\ref{arclength}$)$ 
 has locally bounded second derivative  on $v_f(G)$.
Let $\Hau^1(f(H)) =0$, where $H:=[0,1]\setminus G$.
Suppose that $0 < K < \infty$ and 
$\p$ is an $(f,0,K)$-partition of $G$ such that  $\sum_{I \in \p}\sqrt{V(f,I)} < \infty$. 
Then there exists a homeomorphism $h : [0,1] \to [0,1]$
such that $f \circ h$ is   $D^{2,\infty}$   and  $(f \circ h)'(x)\neq 0$
for each $x \in h^{-1}(G)$. Moreover, if even
\begin{equation}\label{c2}
 F \text{ is}\  C^2 \text{ on}\  v_f(G),
 \end{equation}
then  $f \circ h$ is $C^2$.

Further, in both cases, if $f$ is nonconstant on any interval, then $\lambda(h^{-1}(H)) = 0$.
\end{lemma}

\begin{proof}
Let $U$ be the maximal open set on which $f$ is locally constant. Set $v^*(x) := v_f(x) + \lambda([0,x] \cap U),\ x \in [0,1]$. It is clear that $v^*$ is continuous and increasing. Put $d_1 := v^*(1)$ and $\xi := (v^*)^{-1}$; clearly
 $\xi: [0,d_1] \to [0,1]$ is an increasing homeomorphism. Denote $f_1: = f \circ \xi$, $G_1 := \xi^{-1}(G)$ and $\p_1 := \{\xi^{-1}(I): I \in \p\}$. 
Clearly
\begin{equation}\label{kara}
\sum_{J \in \p_1} \sqrt{\lambda(J)} = \sum_{I \in \p}\sqrt{V(f,I)} < \infty.
\end{equation}
If $(\gamma, \delta)$ is a component of $G_1$, then
   $f_1|_{[\gamma,\delta]}$ is clearly an arc-length parametrization of $f|_{[\xi(\gamma),\xi(\delta)]}$.
    So, using our assumptions  and Lemma \ref{alp}, we obtain that $f_1$ has locally bounded second derivative on $G_1$, and $f_1$ is $C^2$ on $G_1$,
     if \eqref{c2} holds. Further, using  Lemma \ref{alp}, Lemma \ref{arclenlem}, and the properties of $\p$,
      we obtain that
    \begin{equation}\label{I}
    \|f_1'(x)\|=1 \ \ \text{for}\ \  x \in G_1,\ \ \text{and}
    \end{equation}
\begin{equation}\label{LVK}
  \sup _{s \in \inte(J)} \|f_1''(s)\| \cdot \lambda(J)  \leq K\ \ \ \mbox{ for each}\ \ \ J \in
  \p_1.
\end{equation}
Since \eqref{kara} holds, it is easy to see that to each $J \in \p_1$ we can assign  a number
 $\eta_J > 0$ such that
\begin{equation}\label{eta1}
\sum_{J \in \p_1} \sqrt{\lambda(J)/ \eta_J} < \infty,\ \ \ \sqrt{\lambda(J)/ \eta_J} <
1\  (J \in \p_1)\ \ \text{and}
\end{equation}
\begin{equation}\label{eta2}
\{J\in \p_1:\ \eta_J > \ep\}\ \ \text{is finite for each}\ \ \ep>0.
\end{equation}
Denote $d_J := \sqrt{\lambda(J)/ \eta_J}$. For each point $x$ which is an
endpoint of a member of $\p_1$, choose $c(x) \geq 0$ such that
\begin{equation}\label{cex1}
c(x)= 0 \ \ \ \text{whenever}\ \ \ x \notin G_1\ \ \ \text{and}
\end{equation}
\begin{equation}\label{cex2}
0 < c(x) \leq 10^{-10} \cdot (d_J)^2 \eta_J\ \ \text{if}\ \ x\in G_1 \ \ \text{is
an endpoint of some}\ \   J \in \p_1.
\end{equation}
Since $\sum_{J \in \p_1} d_J < \infty$, by Lemma \ref{homint}   we can
 choose   $0 < d_2 < \infty$ and an increasing homeomorphism
 $\Psi: [0,d_2] \to [0,d_1]$ such that $\Psi^{-1}$ is absolutely continuous and  $\lambda(\Psi^{-1}(J)) = d_J $
 for each $J \in \p_1$.  Set  $\p_2 := \{\Psi^{-1}(J):\ J \in \p_1\}$,
 $G_2 := \Psi^{-1}(G_1)$ and $H_2 := [0,d_2] \setminus G_2$.

 Now consider an
 interval $I = [u,v] \in \p_2$. Let $J:= \Psi(I)$, $\eta :=
 \eta_J$, $d:= d_J$, $V := \lambda(J)$, $c_l(I) = c_l := c(\Psi(u))$, and
 $c_r(I) = c_r := c(\Psi(v))$.

Since  $\max(c_l, c_r) \leq 10^{-10} \cdot d^2 \eta$ by~\eqref{cex2} and $d<1$ by \eqref{eta1}, we can choose 
by Lemma \ref{l2}   an increasing $C^2$ homeomorphism
 $\vf_I : I \to J$ such that
\begin{equation}\label{prder}
 \vf_I'(u)=c_l,\  \vf_I'(v) = c_r,\ 0 < \vf_I'(x) \leq  19\cdot \sqrt {\eta V}\quad
 (x \in \inte(I)),
\end{equation}
\begin{equation}\label{drder}
 \vf_I''(u)=\vf_I''(v)=0 \ \ \text{and}\ \
|\vf_I''(x)| \leq  19\cdot \eta\quad (x \in I).
\end{equation}
Thus, setting $\vf(x) := \vf_I(x)$ if $x \in I \in \p_2$
and $\vf(x) := \Psi(x)$ if $x \in [0,d_2]\setminus \bigcup \p_2$,
we easily see that $\vf: [0,d_2] \to [0,d_1]$ is an increasing homeomorphism.
Further, the definition of $c_l(I)$, $c_r(I)$,~\eqref{prder},~\eqref{drder}, and \eqref{cex2} 
easily imply that $\vf$ is $C^2$  on $G_2$ and $\vf'>0$ on $G_2$.
 Since  $(f_1\circ \vf)''(x) =  f_1''(\vf(x))\cdot
(\vf'(x))^2 + f_1'(\vf(x)) \cdot \vf''(x)$ for $x \in G_2$, we easily obtain that
  $f_1 \circ \vf$ is  locally  $D^{2,\infty}$ on $G_2$, and  $f_1 \circ \vf$ is even $C^2$ on $G_2$,
 if  \eqref{c2} holds. Using \eqref{I}, we obtain $(f_1 \circ \vf)'(x) \neq 0$
for $x \in G_2$. 

Now consider an  interval $I = \Psi^{-1}(J) \in \p_2$ and~$x \in  I \cap G_2$. 
Then
\begin{equation}\label{vzor}
 (f_1\circ \vf)''(x) = (f_1 \circ \vf_I)''(x) = f_1''(\vf_I (x))\cdot
(\vf'_I(x))^2 + f_1'(\vf_I(x)) \cdot \vf''_I(x).
\end{equation}   
Consequently, (recall that $\eta = \eta_J$ and $V= \lambda(J)$)
by ~\eqref{I},~\eqref{LVK},~\eqref{prder},~\eqref{drder}, and~\eqref{eta1} we have
\begin{equation}\label{drdersl}
\|(f_1\circ \vf)''(x)\| \leq \frac{K}{V}\cdot 19^2\cdot \eta_JV +
19\cdot \eta_J \leq \eta_J(19^2\cdot K + 19) \ \ \ \text{and}
\end{equation}
\begin{equation}\label{prdersl}
\| (f_1 \circ \vf)'(x)\| = \|f_1'(\vf_I(x))\cdot \vf'_I(x)\| \leq 19
\cdot  \sqrt{\eta_J\lambda(J)}.
\end{equation}
Now we will show that
\begin{equation}\label{suf}
\text{for each}\ \  z \in H_2,\ (f_1 \circ \vf)''(z)=0\ \text{and}\
(f_1 \circ \vf)''\ \text{is continuous at}\ z.
\end{equation}
First, we will prove that $(f_1 \circ \vf)'_+(z)=0$ for each $z \in H_2\setminus \{d_2\}$. 
If $z$ is the left endpoint of some $I \in \p_2$,  we obtain  $(f_1 \circ \vf)'_+(z)=0$    easily from the first equality of~\eqref{prdersl},~\eqref{cex1} 
and the first equality of~\eqref{prder}. If $z$ is not the left endpoint of any $I \in \p_2$, 
 consider an arbitrary $\ep>0$. Using ~\eqref{eta2} and~\eqref{prdersl}, we easily see
that there exists $\delta>0$ such that $\|(f_1 \circ \vf)'(x)\| <\ep$ for each $x \in (z,z+\delta) \cap G_2$. 
Since clearly $\Hau^1((f_1 \circ \vf)(H_2)) =  \Hau^1 (f(H)) =    0$, using the mean value theorem 
and Lemma~\ref{H1L}(iii), we easily obtain that
$f_1 \circ \vf$ is $\ep$-Lipschitz on $[z,z+\delta]$. 
Thus $(f_1 \circ \vf)_+'(z)=0$. Similarly we obtain that
$(f_1 \circ \vf)_-'(z)=0$ for each $z \in H_2 \setminus \{0\}$.
\par
To prove~\eqref{suf}, consider  a point $z \in H_2\setminus \{d_2\}$ and observe that,   
\begin{equation}\label{blizko}
\begin{split}
&\text{for each}\ \ep>0,\ \text{there exists a}\ \delta >0\ \text{such that}\\   
&\|(f_1 \circ \vf)''(x)\|  < \ep\ \ \text{ for each}\ \ x \in (z,z+\delta) \cap G_2.
\end{split}
\end{equation}
Indeed, if $z$ is not the left endpoint of some $I \in \p_2$, then \eqref{blizko}
follows from~\eqref{eta2} and~\eqref{drdersl}.
If $z$ is the left endpoint of some $I \in
\p_2$,  we obtain \eqref{blizko} easily from~\eqref{vzor}, \eqref{I}, \eqref{LVK}, \eqref{cex1},
and the equalities of~\eqref{prder} and~\eqref{drder}.
Now consider an arbitrary $y \in (z,
z+\delta) \cap G_2$; let $w$ be the left endpoint of the component of
$G_2$ which contains $y$. Then $(f_1 \circ \vf)'(w) =0$ and thus the
mean value theorem and \eqref{blizko} imply $\|(f_1 \circ \vf)'(y)\| \leq \ep
|y-w| \leq \ep |y - z|$. Since $(f_1 \circ \vf)'(y) = 0$
for each $y \in (z,z+\delta) \setminus G_2$, we see that
$\lim_{y \to z+} \ (f_1 \circ \vf)'(y)/|y-z| = 0$.
Also using a symmetrical argument, we conclude that $(f_1 \circ
\vf)''(z) =0$ for each $z \in H_2$. The rest
of~\eqref{suf} now  follows from~\eqref{blizko} and its symmetrical version.

We have proved that $(f_1 \circ \vf)''$ is locally bounded on $G_2$. Using also \eqref{suf}, we conclude that
  $(f_1 \circ \vf)''$ is locally bounded on $[0,d_2]$, which implies that $f_1 \circ \vf$ is $D^{2,\infty}$ on  $[0,d_2]$. If \eqref{c2} holds, we have proved that $f_1 \circ \vf$ is $C^2$ on $G_2$; so \eqref{suf} yields
  that  $f_1 \circ \vf$ is  $C^2$ on  $[0,d_2]$.
  
Thus, to finish the proof of the first part of the assertion, it is sufficient to define $h := \xi\circ
\vf \circ \pi$, where $\pi(x) = d_2 x, \ x \in [0,1]$.
\par
To prove the second part of the assertion,
suppose that $f$ is nonconstant on any interval. 
Then $v^*=v_f = \xi^{-1}$; so  
 Lemma \ref{H1L}(iv) implies  $\lambda(\xi^{-1}(H)) = 0$. Since $\Psi^{-1}$ is absolutely continuous 
and $\vf^{-1}(\xi^{-1}(H)) =   \Psi^{-1}(\xi^{-1}(H))$, we have
$\lambda((\xi\circ \vf)^{-1}(H))=0$, and thus also $\lambda(h^{-1}(H))=0$.
\end{proof}

\section{$C^2$-parametrizations}\label{c2sekce}

The following proposition shows that, if $X$  admits an equivalent Fr\'echet smooth norm, the characterization of those $f: [0,1] \to X$,
which allow a $C^2$ parametrization with nonzero  derivatives of the first order, is quite simple.

\begin{proposition}\label{simplepropo}
Let $X$ be a Banach space which admits an equivalent Fr\'echet smooth norm, let
$f:[0,1]\to X$ be continuous. Then the following are equivalent.
\begin{enumerate}
\item 
There exists an increasing homeomorphism $h$
of~$[0,1]$ onto itself such that $f\circ h$ is $C^2$ with $(f\circ h)'(x)\neq0$
for all $x\in[0,1]$. 
\item $f$ is BV, is not constant on any interval, and $F := f \circ v_f^{-1}$ is $C^2$.
\end{enumerate}
\end{proposition}

\begin{proof} 
Let $h$ be as in (i) and  $g :=f\circ h$. It is easy to see that $f$ is BV, is not constant on any interval,
 and $F:= f \circ v_f^{-1} = g \circ v_g^{-1}$. So, (ii) follows from Lemma \ref{vlastF}. 
Now suppose that (ii) holds. Set $\ell:= v_f(1)$ and $h(t):=v_f^{-1}(\ell\cdot t)$, $t\in[0,1]$.
Using Lemma~\ref{arclenlem}, we see that $h$  has the  property from~(i).
\end{proof}

\begin{definition}\label{Df}
For a continuous $f:[0,1]\to X$, denote by {\em $D_f$}  the set of all points $x\in[0,1]$ for which there is no
interval $(c,d)\subset[0,1]$ containing $x$ such
that $f|_{[c,d]}$ has a $C^2$ arc-length parametrization.
\end{definition}
The set $D_f$ is obviously closed and $\{0,1\}\subset D_f$. Further, if $h$ is an increasing homeomorphism of $[0,1]$ onto itself, then
clearly
\begin{equation}\label{Dfgeom}
D_{f \circ h} = h^{-1} (D_f).
\end{equation}

\begin{lemma}\label{Dfvanishlem}
Let $X$ be a Banach space which admits an equivalent Fr\'echet smooth norm,
let $g:[0,1]\to X$ be $C^2$,
and let $x\in D_g$. Then either $g'(x)=0$ or $x\in\{0,1\}$.
\end{lemma}

\begin{proof}
Suppose that $x\in(0,1)\cap D_g$. For a contradiction, assume  $g'(x)\neq0$.
Then there exists $\delta>0$ such that $g'(y) \neq 0$ for each
$y\in[x-\delta,x+\delta]$.
Then Lemma \ref{vlastF} gives that  $g|_{[x-\delta,x+\delta]}$
has a $C^2$ arc-length parametrization, which contradicts $x\in D_g$.
\end{proof}

The following proposition states the main necessary conditions for the existence of a $C^2$ parametrization.

\begin{proposition}\label{nutne}
Let $X$ be a Banach space which admits an equivalent Fr\'echet smooth norm and suppose that a nonconstant continuous 
 $f: [0,1]\to X$ admits a $C^2$ parametrization. Set $G:= (0,1)\setminus D_f$. Then $f$ is BV and the following conditions hold.
 \begin{enumerate}
 
 \item
 $W^{\delta}(f,G)< \infty$ for each $\delta>0$.
 \item
$\sum_{I\in\mcI} \sqrt{V(f,I)} < \infty$, where $\mcI$ is the family   of all components of $G$.
 \item 
 $\int_{v_f(G)}\sqrt{\|F''\|} < \infty$, where $F:=\mcA_f$
(see Definition~\ref{arclength}).
 \item
 $\Hau^1(f(D_f))=0$.
 \end{enumerate}
\end{proposition}
\begin{proof}
Obviously, there exists an increasing  homeomorphism $h$ of $[0,1]$ onto itself
such that $f\circ h$ is a $C^2$ function. Since $f\circ h$ is BV, clearly $f$ is BV as well.
Put $g:=f\circ h$ and $G^* = h^{-1}(G)$. By ~\eqref{Dfgeom}, $D_g = h^{-1}(D_f)$ and $G^* = [0,1] \setminus D_g$.
 Lemma~\ref{Dfvanishlem} implies that $g'(x)=0$ for each $x \in D_g \cap (0,1) = (0,1) \setminus G^*$. 
  Thus $g$ fulfils the assumptions of  Lemma~\ref{Wfin2} (with $f:=g$ and $G:=G^*$). Therefore, for each
   $\delta>0$ we have $W^{\delta}(g,G^*) <\infty$, and so also $W^{\delta}(f,G)< \infty$ by Remark \ref{Wdeltainv}.
  So (i) is proved; 
  (ii) and (iii) follow immediately by  Remark \ref{komponenty} and Lemma \ref{nove}(ii).

 Finally, Lemma~\ref{Sard} implies  $\Hau^1(g(D_g))=\Hau^1(g(D_g \cap (0,1)))=0$. Since $g(D_g) = f(D_f)$, we have proved (iv).

\end{proof}

The main result of the present section is the following theorem which solves the $C^2$ parametrization problem
 for curves in Banach spaces admitting a Fr\' echet smooth norm. 
(Observe
that condition (i) is clearly equivalent to the existence of a $C^2$ parametrization of $f$, and implies that $f$ is  BV).

\begin{theorem}\label{charC2}
Let $X$ be a Banach space which admits an equivalent  Fr\'echet smooth norm.
Let a BV continuous nonconstant $f:[0,1] \to X$ 
be given, and put $G:=[0,1]\setminus D_f$.
Then the following are equivalent.
\begin{enumerate}
\item There exists an increasing homeomorphism $h$ of $[0,1]$ onto itself
such that $f\circ h$ is a $C^2$ function.
\item There exists an increasing homeomorphism $\vf$ of $[0,1]$ onto itself such that
$f\circ\vf$ is $C^2$ and $(f \circ\vf)'(x)\neq 0$ for each $x \in\vf^{-1}(G)$.
\item $\Hau^1(f(D_f))=0$ and $W^{\delta}(f,G)< \infty$ for each $\delta>0$.
\item $\Hau^1(f(D_f))=0$ and $W^{\delta}(f,G)< \infty$ for some $\delta>0$.

\item $\Hau^1(f(D_f))=0$ and $\sum_{I \in \p} \sqrt{V(f,I)} < \infty$
whenever $\p$ is an $(f,\delta,\infty)$-partition of $G$ with $\delta>0$.

\item $\Hau^1(f(D_f))=0$ and $\sum_{I \in \p^*} \sqrt{V(f,I)} < \infty$
for some $(f,0,K^*)$-partition $\p^*$ of $G$ with $0 < K^*<\infty$.

\item  There exists  a family $\s$ of pairwise non-overlapping compact intervals
such that $\inte(J) \subset G$ for each $J \in \s$ and
\begin{equation}\label{tri4}
\sum_{J\in\s} V(f,J)=V(f,[0,1]),\ \
 \sum_{J\in\s} \sqrt{V(f,J)}<\infty,\ \ \sum_{J\in\s}V(f,J)\cdot\sqrt{S_J}<\infty,
\end{equation}
where $S_J=\sup_{y\in v_f(\inte(J))}\|F''(y)\|$ and $F:=\mcA_f$
(see Definition~\ref{arclength}).
\end{enumerate}
If $f$ is nonconstant on any interval, then in~$(ii)$ we can also assert $\lambda(\vf^{-1}(G))=1$.
\end{theorem}

\begin{proof}
 It is sufficient to prove the implications (ii)$\implies$(i)$\implies$(iii)$\implies$(iv)$\implies$(vi)$\implies$
 (vii)$\implies$(iii),\ \ (vi)$\implies$(ii), and (iii)$\implies$(v)$\implies$(vi).

 Note that (since $f$ is nonconstant) both $\Hau^1(f(D_f))=0$ and the equality of~\eqref{tri4}  imply $G \neq \emptyset$
 by Lemma \ref{H1L}(i),(ii).
 
The implication~(ii)$\implies$(i) is trivial and
(i)$\implies$(iii) was proved in Proposition \ref{nutne}.

Obviously, (iii)$\implies$(iv). To prove (iv)$\implies$(vi), let $\delta>0$ be as in condition~(iv), and put $ K^*:=\delta$.
By Lemma~\ref{konstrpartlem} we can choose an $(f,\delta, K^*)$-partition $\p^*$
of $G$. Since $\p^*$ is clearly an $[f,\delta, \infty]$-partition of $G$,  we obtain $\sum_{I\in\p^*} \sqrt{V(f,I)}<\infty$ by Lemma~\ref{nove}(i). Thus we obtain (vi), since  $\p^*$ is clearly an $(f, 0, K^*)$-partition of $G$.

To prove (vi)$\implies$(vii) suppose that~(vi) holds  and $\p^*$ is given. Put $\s := \p^*$.
Then the equality of~\eqref{tri4} holds by Lemma~\ref{H1L}(i) (applied to $G:= \bigcup_{P \in \s}\inte(P)$). 
Since $V(f,J)\cdot \sqrt{S_J} \leq \sqrt{ K^*}\cdot \sqrt{V(f,J)}$ for each $J \in \s$, also both inequalities of~\eqref{tri4} hold.

To prove (vii)$\implies$(iii), let $\s$ be as in (vii). Using Lemma~\ref{H1L}(ii) (with $[a,b] := [0,1]$ and
 $G:= \bigcup_{P \in \s}  \inte{P}$), we obtain   $\Hau^1(f(D_f))\leq  \Hau^1(f([0,1] \setminus  \bigcup_{P \in \s}\inte(P)))=0$. The second condition of (iii) follows immediately from Lemma~\ref{vt2}. 
  
If (vi) holds, then Lemma~\ref{hlavkon2} (with $K:= K^*$ and $H:=D_f$) implies (ii) (and also $\lambda(\vf^{-1}(G))=1$ if $f$
 is nonconstant on any interval).

The implication (iii)$\implies$(v) follows immediately from  Lemma \ref{nove}(i) and Remark \ref{kulhr}. 
 To prove (v)$\implies$(vi), choose by 
 Lemma~\ref{konstrpartlem} an $(f, 1, 1)$-partition  $\p$ of $G$. By (v), $\sum_{I \in \p} \sqrt{V(f,I)} < \infty$; so
  (vi) holds with $\p^* := \p$ and $K^*:=1$.

\end{proof}

\begin{remark}\label{koment}
\begin{enumerate}
\item[(i)]
The assumption that $X$ admits a Fr\' echet smooth norm was used only in the proof of (i)$\implies$(iii). So, any
 of conditions (iii)-(vii) implies (i) and (ii) in an arbitrary $X$.
\item[(ii)]
The conditions (v) and (vi) give an  ``algorithmic'' way how to decide whether (i) holds:

Decide whether $\Hau^1(f(D_f))=0$. If it holds, then choose an
 $(f,\delta,K)$-partition $\p$ of $G := [0,1] \setminus D_f$ with $\delta>0$ and $K< \infty$ (such a partition
  exists by Lemma \ref{konstrpartlem}) and   decide whether      
$\sum_{I \in \p} \sqrt{V(f,I)} < \infty$.

\item[(iii)]
The condition (vii)  needs no auxiliary notions for its formulation. On the other hand, in concrete situations, it is not suitable for a proof that (i) does not hold.
\end{enumerate}
\end{remark}

\begin{remark}\label{nenu3}
Let $X$ be a Banach space which admits an equivalent Fr\'echet smooth norm.
Let $f:[0,1] \to X$ be BV continuous  and 
$G:= [0,1] \setminus D_f$. Then the following are equivalent.
{\it
\begin{enumerate}
\item[(a)] There exists an increasing homeomorphism $\vf$ of $[0,1]$ onto itself such that
$f\circ \vf$ is $C^2$ and  
$(f \circ \vf)'(x)\neq 0$ almost everywhere.
\item [(b)] $f$ is nonconstant on any interval and any of conditions $(iii)-(vii)$ holds.
\end{enumerate}
}
It follows immediately from Theorem~\ref{charC2} 
and the simple observation that~(a) implies that $f$ is nonconstant on any interval. 
\end{remark}

\begin{proposition}\label{monekvivprop}
Let $X$ be a Banach space which admits an equivalent Fr\'echet smooth norm.
Assume that $f:[0,1]\to X$ is continuous, BV and nonconstant on any interval. Let
 $F := f \circ v_f^{-1}$ and $\ell:=v_f(1)$.
Suppose that $F''$ is continuous on $(0,\ell)$   and
there exists $\delta>0$ such that  
$\|F''\|$ is monotone on $(0,\delta)$
and on $(\ell-\delta,\ell)$.
\par
Then the following are equivalent.
\begin{enumerate}
\item
$\int^{\ell}_0 \sqrt{\|F''(t)\|}\,dt<\infty$,
\item
There exists an increasing homeomorphism $h$ of $[0,1]$ onto
itself such that $f\circ h$ is $C^2$.
\end{enumerate}
\end{proposition}

\begin{proof}
The implication (ii)$\implies$(i) follows from Proposition~\ref{nutne}.
 
Now suppose that (i) holds. Choose 
points $(x_i)_{i \in \Z}$ in $(0,\ell)$ such that $\p := \{[x_i,x_{i+1}]:\ i \in \Z\}$ 
is a generalized partition of $(0,\ell)$ and 
$x_{i+1} - x_i\sim  |i|^{-4},\ i \to \pm \infty$. Set $t_i := v_f^{-1}(x_i)$ and 
 $\s := \{v_f^{-1}(I):\ I \in \p\} = \{[t_i,t_{i+1}]:\ i \in \Z\}$.
To prove (ii), it is sufficient to show that ~\eqref{tri4}
from Theorem~\ref{charC2} holds. The first two parts of ~\eqref{tri4} clearly hold; thus it is sufficient to verify that
\begin{equation}\label{treti}
\sum_{J\in\s}V(f,J)\cdot\sqrt{S_J}<\infty,\ \ \ 
\text{where}\ \ \  S_J=\sup_{y\in v_f(\inte(J))}\|F''(y)\|.
\end{equation}
Choose $p\in \N$ such that $x_p > \ell - \delta$ and $x_{-p} < \delta$.  First suppose  that 
 $\|F''\|$ is nondecreasing  on $(\ell-\delta,\ell)$. Then, for each $i \geq p$, 
\[\int^{x_{i+2}}_{x_{i+1}} \sqrt{\|F''(t)\|}\,dt \geq (x_{i+2}-x_{i+1})  \sqrt{\|F''(x_{i+1})\|}
 \geq V(f,[t_{i+1},t_{i+2}]) \sqrt{S_{[t_i,t_{i+1}]}}.\] 
Consequently $\sum_{i=p}^\infty  V(f,[t_{i+1},t_{i+2}]) \sqrt{S_{[t_i,t_{i+1}]}} <\infty$, and since
 $ V(f,[t_{i+1},t_{i+2}]) \sim (i+1)^{-4} \sim i^{-4} \sim V(f,[t_{i},t_{i+1}]),\  i \to \infty$, we obtain
\begin{equation}\label{napravo}
\sum_{i=p}^\infty  V(f,[t_{i},t_{i+1}]) \sqrt{S_{[t_i,t_{i+1}]}} <\infty.
\end{equation}
If  $\|F''\|$ is nonincreasing  on $(\ell-\delta,\ell)$, then  $\|F''\|$ is bounded  on $(\ell-\delta,\ell)$, and thus
 \eqref{napravo} clearly also holds. By a quite symmetrical way we obtain 
\[\sum_{i=p}^\infty  V(f,[t_{-(i+1)},t_{-i}]) \sqrt{S_{[t_{-(i+1)},t_{-i} ]}} <\infty,\]
and therefore 
 \eqref{treti} holds.
\end{proof}

Example~\ref{odm} shows that the assumption on the monotonicity of $\|F''\|$
cannot be omitted in Proposition~\ref{monekvivprop}.
\par
Now we will state a simple consequence of  our results, which might
  be interesting from the point of view of differential geometry. 

\begin{proposition}\label{klas}
Let $\ell \in (0,\infty)$ and $f: (0,\ell) \to \R^n$ be a  $C^2$ curve parametrized by
the arc-length. Let $\kappa_1(s)=\|f''(s)\|$, $s \in (0,\ell) $,
be the first curvature  of $f$. Let $0 <\delta < \infty$ be given.
\par
Then the following are equivalent:
\begin{enumerate}
\item $f$ admits a  $C^2$ parametrization $g: (0,1) \to \R^n$ such that  $g'(x) \neq 0$, $x \in (0,1)$, and $g''$ is bounded.
\item  $V_f^{\delta}:= \sup \{ \sum_{k=1}^n \sqrt{\lambda(I_k)}\} < \infty$, \ \ \ the supremum being taken over all
 systems $(I_k)_{k=1}^n$ of pairwise  nonoverlapping compact subintervals of $(0,\ell)$ such that
  $\lambda(I_k) \max_{s \in I_k} \kappa_1(s) \geq \delta$,\ $k=1,\dots,n$.
\end{enumerate}
If (i) holds, then
\begin{enumerate}
\item[(iii)]  $\int_{0}^{\ell} \sqrt{\kappa_1(s)}\, ds < \infty$.
\end{enumerate}

In the case when $\kappa_1$ is monotone on some right neigbourhood of~$0$ 
and on some left neighbourhood
of~$\ell$, all conditions  (i), (ii) and (iii) are pairwise equivalent.
\end{proposition}

\begin{proof}
Since $f$ is parametrized by
the arc-length, it is
 $1$-Lipschitz; therefore there exists
 a  $1$-Lipschitz  extension $\hat f:[0,\ell] \to X$ of $f$. 
 
 Now consider an arbitrary parametrization $p : [0,1] \to X$ of $\hat f$. Then clearly $\mcA_p = \hat f$. So, using
 only definitions of $V_f^{\delta}$ and $W^{\delta}(p,(0,1))$, we easily obtain
 \begin{equation}\label{srVW}
  V_f^{\delta} \leq W^{\delta}(p,(0,1)) \leq V_f^{\delta} + 2 \sqrt{V(p,[0,1]} = V_f^{\delta} + 2 \sqrt{\ell}.
  \end{equation}
  
  Now suppose that (i) holds, and let $h:(0,1) \to (0,\ell)$ be an increasing homeomorphism such that
   $g = f \circ h$. Let $\hat h: [0,1] \to [0,\ell]$ be the homeomorphism which extends $h$, and set $p:= 
    \hat f \circ \hat h$. Then $p$ is clearly continuous and BV. Since $\mcA_p = \hat f$, we can clearly use
     Lemma \ref{Wfin2} (with $f:= p$ and $G:= (0,1)$) and obtain $W^{\delta}(p, (0,1)) < \infty$. Now
      \eqref{srVW} implies (ii).   Moreover, using Lemma \ref{nove}(ii) (with $f:= p$ and $G:= (0,1)$), we obtain (iii). 
  
  If (ii) holds,  choose  a parametrization $p : [0,1] \to X$ of $\hat f$. Since $D_p = \{0,1\}$, and 
  $W^{\delta}(p,(0,1)) < \infty$ by \eqref{srVW}, we obtain that $p$ admits a $C^2$ parametrization
  $\hat g$ such that $(\hat g)'(x) \neq 0$ for each $x \in (0,1)$ by
  Theorem \ref{charC2} (the equivalence of (ii) and (iv)). Setting $g:= \hat g|_{(0,1)}$, we obtain (i).

Finally, assume that  (iii) holds and that $\kappa_1$ is monotone on $(0,\sigma)$ and $(\ell-\sigma,\ell)$
for some $\sigma>0$. Let  $p : [0,1] \to X$ be a parametrization of $\hat f$. Using Proposition
 \ref{monekvivprop} (with $f:= p$)  and Theorem \ref{charC2} (the equivalence of (i) and (ii)), we easily obtain 
  (i).
\end{proof}
\begin{remark}\label{dalsip}
The proof shows that (i) is equivalent to the condition
\begin{enumerate}
\item[(iv)] There exists a $C^2$ smooth $q: [0,1] \to \R^n$ such that $q'(x) \neq 0$, $x\in (0,1)$, and
 $q|_{(0,1)}$ is a parametrization of $f$.
\end{enumerate}
\end{remark}
\section{Parametrizations with a bounded second derivative}\label{omezddsekce}
Both the results and proofs in the $D^{2,\infty}$ problem  are quite analogous to those in the
 $C^2$ problem. 
 The only two differences are that we now work with $[f,\delta,K]$-partitions 
 instead of $(f,\delta,K)$-partitions, and with the ``singular set'' $\wtilde D_f$ instead of
 $D_f$.
 
 Making these changes (and several other obvious small changes) in the proofs of
 Proposition \ref{simplepropo}, Lemma
  \ref{Dfvanishlem}, Proposition \ref{nutne}, Theorem \ref{charC2} and Proposition \ref{monekvivprop}, we
   obtain
  the proofs of  analogous 
 Proposition \ref{simplepropo2}, Lemma
  \ref{Dfvanishlem2}, Proposition \ref{nutne2}, Theorem \ref{charD2} and Proposition \ref{monekvivprop2}
   below.
    So, we omit the proofs of these results.
 
 The last result of the present section (Proposition \ref{findif}) describes one situation when the  $D^{2,\infty}$ problem
  is equivalent to the $C^2$ problem. 
 
\begin{proposition}\label{simplepropo2}
Let $X$ be a Banach space which admits an equivalent Fr\'echet smooth norm, let
$f:[0,1]\to X$ be continuous. Then the following are equivalent.
\begin{enumerate}
\item 
There exists an increasing homeomorphism $h$
of~$[0,1]$ onto itself such that $f\circ h$ is $D^{2,\infty}$ with $(f\circ h)'(x)\neq0$
for all $x\in[0,1]$. 
\item $f$ is BV, is not constant on any interval, and $F := f \circ v_f^{-1}$ is $D^{2,\infty}$.
\end{enumerate}
\end{proposition}

\begin{definition}\label{Df2}
For a continuous $f:[0,1]\to X$, denote by {\em $\wtilde D_f$}  the set of all points $x\in[0,1]$ for which there is no
interval $(c,d)\subset[0,1]$ containing $x$ such
that $f|_{[c,d]}$ has a $D^{2,\infty}$ arc-length parametrization.
\end{definition}
The set $\wtilde D_f$ is obviously closed and $\{0,1\}\subset \wtilde D_f$. Further, if $h$ is an increasing  homeomorphism of $[0,1]$ onto itself, then
clearly
\begin{equation}\label{Dfgeom2}
\wtilde D_{f \circ h} = h^{-1} (\wtilde D_f).
\end{equation}

\begin{lemma}\label{Dfvanishlem2}
Let $X$ be a Banach space which admits an equivalent Fr\'echet smooth norm,
let $g:[0,1]\to X$ be $D^{2,\infty}$,
and let $x\in \wtilde D_g$. Then either $g'(x)=0$ or $x\in\{0,1\}$.
\end{lemma}

\begin{proposition}\label{nutne2}
Let $X$ be a Banach space which admits an equivalent Fr\'echet smooth norm and suppose that a nonconstant continuous
 $f: [0,1]\to X$ admits an equivalent $D^{2,\infty}$ parametrization. Set $G:= (0,1)\setminus \wtilde D_f$. Then $f$ is BV and     the following conditions hold.
 \begin{enumerate}
 
 \item
 $W^{\delta}(f,G)< \infty$ for each $\delta>0$.
 \item
$\sum_{I\in\mcI} \sqrt{V(f,I)} < \infty$, where $\mcI$ is the family   of all components of $G$.
 \item 
 $\int_{v_f(G)}\sqrt{\|F''\|} < \infty$, where $F:=\mcA_f$
(see Definition~\ref{arclength}).
 \item
 $\Hau^1(f(\wtilde D_f))=0$.
 \end{enumerate}
\end{proposition}

\begin{theorem}\label{charD2}
Let $X$ be a Banach space which admits an equivalent Fr\'echet smooth norm.
Let a BV continuous nonconstant $f:[0,1] \to X$ 
be given, and put $G:=[0,1]\setminus \wtilde D_f$.
Then the following are equivalent.
\begin{enumerate}
\item There exists an increasing homeomorphism $h$ of $[0,1]$ onto itself
such that $f\circ h$ is a $D^{2,\infty}$ function.
\item There exists an increasing homeomorphism $\vf$ of $[0,1]$ onto itself such that
$f\circ\vf$ is $D^{2,\infty}$ and $(f \circ\vf)'(x)\neq 0$ for each $x \in\vf^{-1}(G)$.
\item $\Hau^1(f(\wtilde D_f))=0$ and $W^{\delta}(f,G)< \infty$ for each $\delta>0$.
\item $\Hau^1(f(\wtilde D_f))=0$ and $W^{\delta}(f,G)< \infty$ for some $\delta>0$.

\item $\Hau^1(f(\wtilde D_f))=0$ and $\sum_{I \in \p} \sqrt{V(f,I)} < \infty$
whenever $\p$ is an $[f,\delta,\infty]$-partition of $G$ with $\delta>0$.

\item $\Hau^1(f(\wtilde D_f))=0$ and $\sum_{I \in \p^*} \sqrt{V(f,I)} < \infty$
for some $[f,0,K^*]$-partition $\p^*$ of $G$ with $K^*<\infty$.

\item  There exists  a family $\s$ of pairwise non-overlapping compact intervals
such that $\inte(J) \subset G$ for each $J \in \s$ and
\begin{equation}\label{tri42}
\sum_{J\in\s} V(f,J)=V(f,[0,1]),\ \
 \sum_{J\in\s} \sqrt{V(f,J)}<\infty,\ \ \sum_{J\in\s}V(f,J)\cdot\sqrt{S_J}<\infty,
\end{equation}
where $S_J=\sup_{y\in v_f(\inte(J))}\|F''(y)\|$ and $F:=\mcA_f$
(see Definition~\ref{arclength}).
\end{enumerate}
If $f$ is nonconstant on any interval, then in~$(ii)$ we can also assert $\lambda(\vf^{-1}(G))=1$.
\end{theorem}

Let us note that the remarks quite analogous 
 to Remark \ref{koment} and Remark \ref{nenu3} (we only replace $C^2$ by $D^{2,\infty}$, $(f,\delta,K)$ by
 $[f,\delta,K]$, and $D_f$ by $\wtilde D_f$)
 are true.

\begin{proposition}\label{monekvivprop2}
Let $X$ be a Banach space which admits an equivalent Fr\'echet smooth norm.
Assume that $f:[0,1]\to X$ is continuous, BV and nonconstant on any interval. Let
 $F := f \circ v_f^{-1}$ and $\ell:=v_f(1)$.
Suppose that $F''$ is locally bounded on $(0,\ell)$   and
for some $\delta>0$ we have that  
$\|F''\|$ is monotone on $(0,\delta)$
and on $(\ell-\delta,\ell)$.
\par
Then the following are equivalent.
\begin{enumerate}
\item
$\int^{\ell}_0 \sqrt{\|F''(t)\|}\,dt<\infty$,
\item
There exists an increasing homeomorphism $h$ of $[0,1]$ onto
itself such that $f\circ h$ is $D^{2,\infty}$.
\end{enumerate}
\end{proposition}

Our main results easily imply the following proposition, which should be compared with
 Example~\ref{dmd} and  Example~\ref{dmd2} below.

\begin{proposition}\label{findif}
Let $f: [0,1] \to X$ be continuous and BV such that $B:= D_f \setminus \wtilde D_f$ is finite. Then $f$
 admits a $C^2$ parametrization if and only if $f$ admits a $D^{2,\infty}$ parametrization.
\end{proposition}
\begin{proof}
The ``only if'' implication is obvious. So, suppose that $f$ admits a $D^{2,\infty}$ parametrization
and choose a family $\s$ by condition (vii) of Theorem \ref{charD2}. Dividing members of $\s$ by points of
 $B$, we obtain a non-overlapping system $\s^*$ of compact intervals (such that
 $\bigcup_{I \in \s^*} \inte(I) \subset
   \bigcup_{J \in \s} \inte(J)$ and $\bigcup_{I \in \s} \inte(I) \setminus
   \bigcup_{J \in \s^*} \inte(J) =  B \cap  \bigcup_{I \in \s} \inte(I)$). Since $\s^*$ clearly shows that the condition (vii) of Theorem \ref{charC2}
   holds, $f$ admits a $C^2$ parametrization by Theorem \ref{charC2}.
\end{proof}

\section{Examples}\label{examples}

\begin{example}\label{dmd}
There exists a  function  $f: [0,1] \to \R^2$ which admits a $D^{2,\infty}$ parametrization, does not 
 admits a $C^2$ parametrization and $ D_f \setminus \wtilde D_f = \{1/n:\ n= 2,3,\dots\}$. 
\end{example}
 \begin{proof}
 Set $C:= \{1/n:\ n=2,3,\dots\} \cup \{0, 1\}$. By Lemma \ref{bru} there exists a real function $\vf$ on $[0,1]$ which 
 has a bounded derivative on $[0,1]$, and $C$ is the set of discontinuity of~$\vf'$.
 Now identify $\R^2$ with the complex plane and set $f(x) := \int_0^x e^{i \vf(t)}\ dt$. Then $f$ is Lipschitz,
    $\|f'(x)\|=1$, $x \in [0,1]$, and thus $f$ is parametrized by the arc-length by \eqref{varder}. Since
     $f''(t) = ie^{i\vf(t)} \vf'(t)$ for $t \in (0,1)$, we easily see that $\wtilde D_f = \{0,1\}$ and 
      $D_f = C$.
       Since  the condition
        (vii) of Theorem \ref{charD2} clearly holds with $\s := \{[0,1]\}$,  $f$ admits a $D^{2,\infty}$ parametrization. On the other hand, $\sum_{n\in \N} \sqrt{V(f,[1/(n+1),1/n])} = \sum_{n\in \N} \sqrt{1/(n(n+1))} = \infty$. So, the condition (ii) of Proposition \ref{nutne}  does not hold, and consequently
  $f$ does not admit a $C^2$ parametrization.       
 \end{proof}

\begin{example}\label{dmd2}
There exists a  function  $f: [0,1] \to \R^2$ which 
 admits a $C^2$ parametrization and $ D_f \setminus \wtilde D_f $ is uncountable.
 \end{example} 
 \begin{proof} 
 Let $C \subset [0,1]$ be a closed set which is constructed by the same way as the Cantor ternary set, 
 with the only difference that in the $n$-th step of the construction  from any of $2^{n-1}$ closed intervals $I$ we delete
  a concentric open interval of the length  $(3/5) \lambda(I)$. So, $H := (0,1) \setminus C$ has, for each $n \in \N$, 
  $2^{n-1}$ components of length $(3/5)(1/5)^{n-1}$ (and no other component). So, 
   $\lambda(C) = 1 - \sum_{n=1}^{\infty}  (3/5)(2/5)^{n-1} =0$. By Lemma~\ref{bru} there exists a real function $\vf$ on $[0,1]$ which has a bounded derivative on $[0,1]$, and $C$ is the set of discontinuity of $\vf'$.
   Setting  $f(x) := \int_0^x e^{i \vf(t)}\ dt$, we obtain as in Example~\ref{dmd} that $f$ is parametrized by the arc-length,  $\wtilde D_f = \{0,1\}$ and $D_f = C$. 
   To prove that $f$ admits a $C^2$ parametrization, we will
    verify condition (vii) of Theorem~\ref{charC2}. 
  We will show that~\eqref{tri4} holds, if we define $\s$ as the system of the closures of all components of $H$. 
First, \[ \sum_{J \in \s} V(f,J) = \sum_{J \in \s} \lambda(J) =1 = V(f, [0,1]).\]
    Further,
    \[ \sum_{J \in \s} \sqrt{V(f,J)} = \sum_{n=1}^{\infty}  2^{n-1} \sqrt{(3/5)(1/5)^{n-1}} < \infty.\]
    Finally, the last part of \eqref{tri4} also holds, since $S_J =1$ for each $J \in \s$.
 \end{proof}

\begin{example}
For $s>0$, consider the spiral $f: [0,1] \to \R^2$ defined by $f(0)=0$ and
\[ f(t) = (x(t), y(t)) = (t^s \cos(1/t), t^s \sin (1/t)),\ \ \ 0< t \leq 1.\]
Identifying $\R^2$ with the complex plane, we have, for $t \in (0,1)$,
\[  f(t) = t^s e^{i/t}\ \ \ \text{and}\ \ \ f'(t) = e^{i/t} t^{s-2} (st - i).\]
 Consequently,  $\|f'(t)\| \sim t^{s-2},\ t \to 0+$.
  So, using \eqref{varder} on intervals $[\delta,1]$,
  we easily obtain that $f$ is BV if and only if $s>1$. 
Using the well-known formula for the oriented curvature $k(t)$ (see e.g. \cite[p.~26]{S}) we obtain
\begin{equation}\label{krivost}
 k(t) = \frac{\det (f'(t), f''(t))}{\|f'(t)\|^3} = \frac{s(1-s)t^{2s-4} - t^{2s-6}}{t^{3s-6} (1+s^2t^2)^{3/2}} < 0,\quad t\in (0,1).
\end{equation}
Denote $F:= f \circ v_f^{-1}$ and $\ell:= v_f(1)$. Since $\|F''(v_f(t))\| = |k(t)|$ and $v_f(t)$ is increasing,
  \eqref{krivost} easily implies that there is a $\delta>0$ such that $\|F''\|$ is monotone on
   $(0,\delta)$ and on $(\ell-\delta,\ell)$. Since $v_f'(t) = \|f'(t)\|$, we have
   \[  \int_0^\ell\ \sqrt {\|F''(y)\|}\ dy =  \int_0^1 \ \sqrt{|k(t)|} \, \|f'(t)\|\ dt.\]
Using \eqref{krivost}, we easily obtain $\sqrt{|k(t)|} \, \|f'(t)\| \sim t^{s/2-2},\ t\to 0+\, $. Consequently, 
 using Proposition \ref{monekvivprop} (resp. Proposition \ref{monekvivprop2}), we obtain that
  $f$ admits a $C^2$ (resp. $D^{2,\infty}$) parametrization if and only if $s>2$.  
\end{example}

In the following example, we need the following well-known fact.

\begin{lemma}\label{dglem} 
Let $k:(0,1)\to\R$ be positive and $C^\infty$.
Then there is a continuous $f:[0,1]\to\R^2$ parametrized
by the arc-length, $C^\infty$ on $(0,1)$, and such that $\|f''(x)\|=k(x)$
for  $x\in(0,1)$.
\end{lemma}
\begin{proof}
By the Fundamental Theorem of the local theory of curves (see e.g.\ \cite[Theorem~2.15]{Kuh}),
there exists  $g:(0,1)\to\R^2$ parametrized by the arc-length, which is
$C^\infty$, and  $\|g''(x)\|=k(x)$ for $x\in(0,1)$.
Since $g$ is $1$-Lipschitz, it has a continuous extension $f$  to 
$[0,1]$, which has all the desired properties.
\end{proof}

\begin{example}\label{odm}
There exists a continuous $f:[0,1]\to\R^2$ which is parametrized by the 
arc-length, is $C^{\infty}$ on $(0,1)$, does not allow a $D^{2,\infty}$-parametrization,
but $\int_0^1 \sqrt{\|f''\|} < \infty$.
\end{example}

\begin{proof}
Let $\p = \{I_n:\ n \in \N\}$ be a generalized partition of $(0,1)$ such that
$\lambda(I_n) =c/n^2$ for some $c>0$. Choose closed intervals $J_n \subset I_n$ with  $\lambda(J_n) = c/n^4$. 
Clearly, we can choose~a positive $C^\infty$ function $k:(0,1)\to (0,\infty)$ such that, for each $n \in \N$, 
we have
$\max_{x\in J_n} k(x)= n^4$ and  $k(x)\leq1$ for  $x \in I_n \setminus J_n$.
Choose an $f$ corresponding to $k$ by Lemma~\ref{dglem}. 
We easily see that $\p$ is an $(f,c,\infty)$-partition of~$(0,1)$,
but $\sum_{I \in \p} \sqrt{V(f,I)} = \sum_{I \in \p} \sqrt{\lambda(I)} = \infty$. So, condition (v) of Theorem~\ref{charD2} does not hold (see Remark \ref{kulhr}), and therefore $f$ does not allow
a $D^{2,\infty}$-parametrization. 
On the other hand, $\int_0^1 \sqrt{\|f''\|} \leq  1 + \sum (c/n^4)\sqrt{n^4} < \infty$.
\end{proof}

\section{The case of real valued functions}

As we noted in Introduction, the case of ``higher order smooth'' parametrizations for
$X=\R$ was settled independently by Laczkovich and Preiss \cite{LP}
and Lebedev~\cite{Leb}. Both papers contain (formally slightly
different) characterizations of those $f: [0,1] \to \R$ which
allow an equivalent $C^n$ ($n \in \N)$ parametrization (or a parametrization with bounded $n$-th derivative).
 Lebedev's results give that  for a continuous $f:[0,1]\to\R$ the following conditions are equivalent:

\begin{enumerate}
\item [(i)]
$f$ admits a $C^n$ parametrization.
\item[(ii)]
$f$ admits a parametrization with bounded $n$-th  derivative.
\item[(iii)]
\[
\lambda(f(K_f))=0\ \text{ and}\ 
\sum_{\alpha \in A} (\omega^f_{\alpha})^{1/n}<\infty, 
\]
where $(I_{\alpha})_{\alpha\in A}$ are all maximal open intervals in $[0,1]$ on which $f$ is
constant or strictly monotone,   $K_f := [0,1]\setminus \bigcup_{\alpha\in A} I_{\alpha}$
is the set of points ``of varying monotonicity'' of $f$ (denoted by $M_f$ in \cite{Leb}),
and $\omega^f_{\alpha}$ is the oscillation of $f$ on $I_{\alpha}$.
\end{enumerate}
\bigskip
Laczkovich and Preiss~\cite{LP} show that the conditions (i) and (ii) are equivalent to
\begin{enumerate}
\item[(iv)]
\[
V_{1/n}(f,K_f) < \infty,
\]
where 
\[ V_{1/n}(f,K_f) : = 
\sup \{\sum_{i=1}^m |f(d_i)-f(c_i)|^{1/n}\},\]
the supremum being taken over all systems $[c_i,d_i]$, $i=1,\dots,m$,
of pairwise non-overlapping subintervals of $[0,1]$ with $c_i, d_i \in K_f$.
\end{enumerate}

We will indicate how, {\it in the case $n=2$}, the equivalence of conditions (i), (ii) and (iii)  follows from the results
of the present article (without using any result  of~\cite{LP} and ~\cite{Leb}).
 First we will show that
\begin{equation}\label{ssr}
  K_f \cup E_f =  \wtilde D_f = D_f \ \ \ \text{for each continuous}\ \ \ f:[0,1] \to \R,
\end{equation}
where $E_f$ is the maximal open subset of $(0,1)$ on which $f$ is locally constant.
To prove $(K_f  \cup E_f) \subset \wtilde D_f$, suppose that $x \in [0,1]\setminus \wtilde D_f$. By Definition \ref{Df2} 
 there is an
interval $(c,d)\subset[0,1]$ containing $x$ such
that $f|_{[c,d]}$ has a $D^{2,\infty}$ arc-length parametrization $f^*$. Since $f^* \in C^1$, we obtain that
 $|f^*|=1$ on the domain of $f^*$ by  Lemma
 \ref{arclenlem}. Thus $f^*$ is affine with the slope $1$ or $-1$. Consequently, $f$ is strictly monotone
  on $(c,d)$, which implies $x \notin  K_f  \cup E_f$. The inclusion $\wtilde D_f  \subset D_f$ is obvious.
   To prove $  D_f  \subset  K_f  \cup E_f$, suppose that $x \in [0,1]\setminus (K_f  \cup E_f) $. Then
    there is an
interval $(c,d)\subset[0,1]$ containing $x$ such that $f$ is strictly monotone on $[c,d]$. Since 
 $V(f,[c,x]) = \pm (f(x)-f(c))$ for $x \in [c,d]$, we easily obtain that each arc-length parametrization
  of  $f|_{[c,d]}$ is affine with the slope $1$ or $-1$. So, $x \notin D_f$.

The implication (i)$\implies$(ii) is trivial. 

Now suppose that (ii) holds for $n=2$. By Proposition \ref{nutne2}(iv) and \eqref{ssr} we obtain
 $\lambda(f(K_f))=0$. 
 Let $(I_{\alpha})_{\alpha\in A}$ be as in (iii) and let $A^*$ be the set
 of those $\alpha \in A$, for which $f$ is strictly monotone on $I_{\alpha}$. By \eqref{ssr} we see
  that  $(I_{\alpha})_{\alpha\in A^*}$ is the system of all components of $(0,1) \setminus \wtilde D_f$.
  So, using Proposition \ref{nutne2}(ii), we obtain $\sum_{\alpha \in A} (\omega^f_{\alpha})^{1/2} =
 \sum_{\alpha \in A^*} (\omega^f_{\alpha})^{1/2} < \infty$. Thus we have proved (iii).
 
 Finally suppose that (iii) holds for $n=2$, and let  $I_{\alpha},\ \alpha \in A,$ be as in (iii). Since $\omega_{\alpha}^f = V(f, \overline{I_{\alpha}})$, we obtain
 $\sum_{\alpha \in A} V(f, \overline{I_{\alpha}}) < \infty$, and therefore $f$ is BV by Lemma~\ref{H1L}(i).
 We will prove that condition~(vii) of Theorem~\ref{charC2} is satisfied. Let $\s$ be the family 
of those $\overline{I_{\alpha}} $ on which~$f$ is not constant. The equality \eqref{ssr} gives that
 $\inte(J) \subset G:= [0,1] \setminus D_f$ for each $J \in \s$.
  Since clearly $\lambda(f(E_f))=0$, we obtain
  $\lambda(f(D_f))= \lambda(f(K_f))=0$ by \eqref{ssr}.
Consequently Lemma~\ref{H1L}(i) implies 
$\sum_{J\in \s} V(f,J)= \sum_{\alpha \in A} V(f,I_{\alpha}) =   V(f,[0,1])$, 
which is the first equality of~\eqref{tri4}. 
The validity of the second equality of \eqref{tri4} immediately follows from (iii).
Since $f$ is strictly monotone on each~$J \in \s$, 
we obtain (as in the proof of \eqref{ssr}) that $F:=\mcA_f$
(see Definition~\ref{arclength}) is affine on $v_f(J)$. 
So, $\s_J = 0$, and the third equality of~\eqref{tri4} holds.
\medskip
 
\begin{remark}\label{LPLeb}
Finally note that, using  Lemma~\ref{H1L}(i),(ii), it is easy to give a direct proof 
that the Laczkovich-Preiss condition (iv)
is equivalent to the Lebedev condition (iii) (for any $n \in \N$). However, since this  proof is 
 not short, we do not present it here.
 \end{remark}

\section*{Acknowledgment}
The research of the first author was supported in part by ISF.
The research of the second author  was supported by
the institutional grant MSM 0021620839 and the grants GA\v CR  201/03/0931   and  GA\v CR 201/06/0198. 
Both authors would also like to thank for hospitality to Erwin Schr\"odinger Institute, Vienna,
where part of the work was done.

\end{document}